\documentclass[a4paper,10pt]{amsart}

\usepackage[all]{xy}
\usepackage{amssymb}
\usepackage[hypertex]{hyperref}
\newcommand{\Spe}{{\rm Spec}}

\newcommand{\tg}{{\rm tg}}
\newcommand{\Dgl}{{D_{\rm gl}}}

\newcommand{\TO}{{\mathcal{T}_\mathcal{O} }}
\newcommand{\T}{{\mathcal{T}}}
\newcommand{\A}{{ \rm Aut }}

\newcommand{\F}{{\mathbb{F}}}

\renewcommand{\O}{{\mathcal{O}}}
\newcommand{\Z}{{\mathbb{Z}}}
\newcommand{\lc}{\left\lceil}
\newcommand{\rc}{\right\rceil}
\newcommand{\lf}{\left\lfloor}
\newcommand{\rf}{\right\rfloor}
 
\newcommand{\Gal}{{ \rm Gal}}
\renewcommand{\mod}{{\;\rm mod}}
\newtheorem{pro}{Proposition}[section]
\newtheorem{lemma}[pro]{Lemma}
\newtheorem{rem}[pro]{Remark}
\newtheorem{cor}[pro]{Corollary}
\newtheorem{theorem}[pro]{Theorem}

\begin{document}
\bibliographystyle{amsplain}

\title{On the tangent space of the deformation functor 
of Curves with Automorphisms}

\author{A. Kontogeorgis}

\begin{abstract} We provide a method to 
compute the  dimension of the tangent space to  the 
global infinitesimal deformation functor of  
a curve together with a subgroup of 
the group of automorphisms. The computational techniques we 
developed are applied to  several examples including Fermat curves, $p$-cyclic covers 
of the affine line and to Lehr-Matignon curves.
\end{abstract}

\email{kontogar@aegean.gr}
\address{
Department of Mathematics, University of the \AE gean, 83200 Karlovassi, Samos,
Greece}
\date{\today}

\maketitle

%opening
The aim of this paper is  the  study  of   equivariant equicharacteristic 
infinitesimal deformations of a curve $X$ of genus $g$, 
admitting a group of automorphisms. This paper is the result of 
my attempt to 
understand the works of J.Bertin - A.M\'ezard \cite{Be-Me} 
and of  G. Cornelissen - F. Kato \cite{CK}. 

%The authors of the above papers, investigate infinitesimal  deformations 
%of a curve  compatible with the action of a subgroup of the group of
%automorphisms. 

Let $X$ be a smooth  algebraic curve, defined over an algebraically closed 
field of characteristic $p\geq 0$. 
The infinitesimal deformations of the  curve $X$, 
without considering compatibility with the group action, 
correspond to  directions on the vector space
$H^1(X,\T_X)$  which constitutes the tangent space to the 
 deformation functor of  the curve $X$ \cite{HarrisModuli}. 
All elements in $H^1(X,\T_X)$ give rise to unobstructed 
deformations, since $X$ is one-dimensional and the second 
cohomology vanishes. 

In the study of deformations together with the action
of a subgroup of the automorphism group, 
a new deformation functor can be defined. The tangent space of this 
functor  is given 
by  Grothendieck's  \cite{GroTo} equivariant cohomology  group $H^1(X,G,\T_X)$,  
   \cite[3.1]{Be-Me}. 
 In this case 
 the wild ramification points
contribute to the dimension of the tangent space of the deformation  functor
and also posed several lifting obstructions, related to the 
theory of 
deformations of Galois representations.

The authors of \cite{Be-Me}, after proving a 
local-global principle, focused on infinitesimal deformations 
in the case $G$ is cyclic of order $p$ and considered liftings
to characteristic zero, while the authors of \cite{CK} considered 
the case of deformations of ordinary curves without putting any other
condition on the 
automorphism group. The ramification groups of  automorphism groups
acting on ordinary curves have a special ramification filtration, 
{\em i.e.}, the $p$-part of every 
ramification group
 is an elementary Abelian group, and this 
makes the computation possible, since  elementary Abelian 
group extensions are   given explicitly in terms of   Artin-Schreier 
extensions.  

In this paper we consider an arbitrary curve $X$ with automorphism 
group $G$. By the theory of 
Galois groups of local fields,
the ramification group at every wild ramified point, 
can break to a sequence of extensions of elementary Abelian 
groups \cite[IV]{SeL}. We will use this decomposition,  together with the
spectral sequence of Lyndon-Hochschild-Serre in order to 
reduce the computation,  to a computation 
involving  elementary Abelian groups.

We are working over an algebraically closed field of positive characterstic and 
for the sake of simplicity we assume that $p\geq 5$.

The dimension of the tangent space of the deformation 
functor, depends on the group structure 
of the extensions that appear in the series decomposition of 
the ramification  groups at wild ramified points. 
We are able to give lower
and upper bounds  of the dimension
of the tangent space of the  deformation functor.
%As far as the author knows the computation 
%in this case,  remains open even in the cyclic case. 
 
In particular, if the decomposition group $G_P$ at a wild ramified point $P$ is the semidirect product of an 
elementary abelian group with a cyclic group, such that there is only a lower jump at the $i$-th position 
in the ramification filtration,  then we are able to 
compute exactly the dimension of the local contribution  
$H^1(G_P,\TO)$ prop. \ref{comp-el-ab} and example 4 on page \pageref{examp4}.

We begin our exposition in section \ref{surveydef} by surveying some of the known 
deformation theory. 
Next we proceed to the most difficult task, namely the  computation of the tangent space of the local 
deformation functor, by employing the low terms sequence stemming from the 
Lyndon-Hochschild-Serre spectral sequence.  

The dimension of equivariant deformations that are locally trivial,
{\em i.e.}, the dimension of $H^1(X/G, \pi_*^G (\T_X))$ is computed in section \ref{global-comp}.
The computational techniques we developed are applied to 
the case of Fermat curves, that are known to have 
large automorphism group, to the case of $p$-covers of $\mathbb{P}^1(k)$ and to the case 
of Lehr-Matignon curves. 
Moreover, we are able to recover the results of Cornelissen-Kato \cite{CK} concerning 
deformations of ordinary  curves. Finally, we try to compare our result with the 
results of R. Pries \cite{Pries:02},\cite{Pries:04} concerning the computation of unobstructed deformations of wild ramified actions on curves.  

{\bf Acknowledgement} The author would like  to thank R. Pries, for her  useful comments and corrections.

%The study of the hull, and the computation of its Krull dimension, 
%involves the study of lifting obstructions
%and this will the subject of a forthcoming paper. 

%
%
%
\section{Some Deformation Theory} \label{surveydef}
%
%
%
%
%  A survey of Deformation Theory
% 
There is nothing original in this section, but 
for the sake of completeness, 
we present some of the tools we will need for our study.
This part is essentially a review of  \cite{Be-Me},\cite{CK}, 
\cite{MazDef}.

Let $k$ be an algebraic closed field  of characteristic $p\geq 0$. 
We consider the category $\mathcal{C}$ of local Artin $k$-algebras 
with residue field $k$. 

Let $X$ be a non-singular  projective curve defined over the field $k$, and let $G$ be a 
 fixed subgroup of the automorphism group of $X$. 
We will denote by  $(X,G)$  the couple of the curve $X$ together with the  
group $G$. 

A deformation of the couple $(X,G)$ over the local Artin ring $A$ is a proper, smooth family
of curves 
\[
\mathcal{X} \rightarrow  \Spe(A)
\]
parametrized by the base scheme $\Spe(A)$, together with a group homomorphism $G\rightarrow \A _A(\mathcal{X})$ such that there is a 
$G$-equivariant isomorphism $\phi$ 
from 
the fibre over the closed 
point of $A$ to the original curve $X$:
\[
\phi: \mathcal{X}\otimes_{\Spe(A)} \Spe(k)\rightarrow X. 
\]
Two deformations $\mathcal{X}_1,\mathcal{X}_2$ are considered 
to be equivalent if there is a $G$-equivariant isomorphism $\psi$, 
making the following diagram  commutative:
\[
\xymatrix{
\mathcal{X}_1 \ar[rr]^{\psi} \ar[dr] & & \mathcal{X}_2 \ar[dl] \\
& \Spe A &
}
\]

The  global deformation functor is defined:
\[
\Dgl: \mathcal{C} \rightarrow \rm{Sets}, 
A \mapsto
\left\{
\mbox{
\begin{tabular}{l}
Equivalence classes \\
of deformations of \\
couples $(X,G)$ over $A$
\end{tabular}
}
\right\}
\]
Let $D$ be a functor such that $D(k)$ is a single element. 
If $k[\epsilon]$ is the ring of dual numbers, then the 
Zariski tangent space $t_D$  of the functor is defined by  
$t_D:=D(k[\epsilon])$. 
If the functor $D$ satisfies the ``Tangent Space Hypothesis", 
{\em i.e.}, when the mapping 
\[
h:D(k[\epsilon]\times_k k[\epsilon]) \rightarrow D(k[\epsilon]) \times D(k[\epsilon])
\]
is an isomorphism,
then the $D(k[\epsilon])$ admits the structure of a $k$-vector space \cite[p.272]{MazDef}. 
The tangent space hypothesis is contained in $(H3)$-Schlessinger's hypothesis that 
hold for all the functors in this paper, since all the functors admit versal deformation rings 
\cite{Sch},\cite[sect. 2]{Be-Me}. 

The tangent space $t_{D_{gl}}:=\Dgl(k[\epsilon])$ of the global deformation 
functor is expressed  in terms of Grothendieck's equivariant 
cohomology, that combines the 
construction of group  cohomology and  sheaf cohomology \cite{GroTo}. 

We recall quickly the definition of equivariant cohomology theory:
We consider the covering map $\pi: X \rightarrow Y=X/G$. 
For every sheaf $F$ on $X$ we denote by $\pi_*^G(F)$ the sheaf 
\[
V\mapsto \Gamma(\pi^{-1}(V), F)^G, \mbox{ where  } V \mbox{ is an open set of } Y
\]
The category of $(G,\mathcal{O}_X)$-modules is the category 
of $\mathcal{O}_X$-modules with an additional  $G$-module structure. 
We can define two left exact functors from  the category of 
$(G,\mathcal{O}_X)$-modules, namely
\[
\pi_*^G \mbox{ and } \Gamma^G(X, \cdot),
\]
where $\Gamma^G(X,F)=\Gamma(X,F)^G$. The derived functors 
$R^q \pi_*^G(X,\cdot)$
of the first 
functor are sheaves of modules on $Y$, and the derived functors of the 
second are groups $H^q(X,G,F)=R^q \Gamma^G (X,F)$. 
%
%% 5 Dec 2001
%
J. Bertin and A. M\'ezard  in \cite{Be-Me} proved the following 
\begin{theorem}
Let $\T_X$ be the tangent sheaf on the curve $X$.
The tangent space $t_{\Dgl}$ to the global deformation functor, 
is given in 
terms of equivariant cohomology as $t_{\Dgl}=H^1(X,G,\T_X)$. 
Moreover the following sequence is exact:
\begin{equation} \label{tang-equiv-gal}
0 \rightarrow H^1(X/G, \pi_*^G(\T_X)) \rightarrow
H^1(X,G,\T_X) \rightarrow 
H^0(X/G,  R^1 \pi_*^G(\T_X)) \rightarrow 0.
\end{equation}
\end{theorem}

For a local ring $k[[t]]$ we define the 
local tangent space $\TO$, as the $k[[t]]$-module of 
$k$-derivations. The module $\TO:=k[[t]] \frac{d}{dt}$, where 
$\delta=\frac{d}{dt}$ is the derivation such that $\delta(t)=1$. 
If $G$ is a subgroup of $\A (k[[t]])$, then $G$ acts on 
$\TO$ in terms of the adjoint representation. 
Moreover there is a bijection $D_\rho(k[\epsilon])\stackrel{\cong}{\rightarrow}
H^1(G,\TO)$ \cite{CK}.

In order to describe the tangent space of the local  deformation space we 
will   compute 
first the space of tangential liftings, {\em i.e.}, the 
space $H^1(G,\TO)$. 

This problem is solved when $G$ is a cyclic group of order $p$, 
by  J.Bertin and A. M\'ezard  \cite{Be-Me} and if the original curve is ordinary 
by G. Cornelissen and F. Kato in \cite{CK}. 
%The authors of the above papers 
%succeeded also in the computation of the lifting obstructions.   

We will apply the classification of groups 
that can appear as Galois groups of local fields 
in order to reduce the problem to elementary Abelian group  case.

\subsection{Splitting the branch locus} \label{split-branch-points}
Let $P$ be a wild ramified point on the special fibre $X$, and let 
$\sigma \in G_j(P)$ where $G_j(P)$ denotes the $j$-ramification group at $P$. Assume that we can deform  the special fibre 
to a deformation $\mathcal{X} \rightarrow A$, where $A$ is a complete 
local discrete valued ring that is a $k$-algebra. Denote  by  
$m_A$ the maximal ideal of $A$ and assume that $A/m_A=k$. Moreover assume   that $\sigma$ acts fibrewise on 
$\mathcal{X}$. We will follow Green-Matignon \cite{Green-Mat}, on 
expressing the expansion 
\[
\sigma(T)-T= f_j(T) u(T),
\]   
where $f_j(T)=\sum_{\nu=0}^j a_i T^i$ ($a_i \in m_A$ for $\nu=0,\ldots,j-1$, $a_j=1$)  is a distinguished Weierstrass polynomial of 
degree $j$ \cite[VII 8. prop. 6]{BourbakiComm} and $u(T)$ is a unit of 
$A[[T]]$. 
The reduction of the polynomial $f_j$ modulo $m_A$
gives the automorphism $\sigma$ on $G_j(P)$ but 
$\sigma$  when lifted on $\mathcal{X}$ has in general more than one fixed 
points, since $f_j(T)$ might be a reducible polynomial. 
If $f_j(T)$, gives rise to only one  horizontal branch divisor
then we say that the corresponding deformation does 
not split the branch locus. 
%This situation is rather special,
%since it puts various discriminant conditions on the 
%coefficients $a_i$ of  the polynomial $f_j$.  
%INERTIA!!!

Moreover,
if we reduce $\mathcal{X} \times_A  \Spe \frac{A}{m_A^2}$ we obtain an 
infinitesimal extension that gives rise to a cohomology 
class in $H^1(G(P),\TO)$ by \cite[prop. 2.3]{CK}.

On the other hand cohomology classes in $H^1(X/G, \pi_*^G(\T_X))$ induce trivial 
deformations on formal neighbourhoods of the branch point $P$ \cite[3.3.1]{Be-Me}
and do not split the branch points.  In the special case of ordinary curves, the 
distinction of deformations that do or do not split the branch points does not 
occur since the polynomials $f_j$ are of degree $1$.

\subsection{Description of the ramification Group}
The finite  groups that appear as Galois groups of 
a local field $k((t))$, where $k$ is algebraically 
closed of characteristic $p$ are known \cite{SeL}. 

Let $L/K$ be a Galois  extension of a local field $K$ with 
Galois group $G$. 
We consider  the ramification filtration of $G$, 
\begin{equation}\label{ram-filt}
G=G_0 \subset G_1 \subseteq  G_2 \subseteq \cdots \subseteq G_n \subset 
G_{n+1} =\{1\}.
\end{equation}
The quotient  $G_0/G_1$ is a cyclic group of order prime to the 
characteristic, $G_1$ is $p$-group and for 
$i\geq 1$ the quotients
$G_i/G_{i+1}$ are  elementary Abelian $p$-groups. 
If a curve is ordinary then by \cite{Nak} the ramification 
filtration is short, {\em i.e.},  $G_2=\{1\}$, and this gives that 
$G_1$ is an elementary Abelian group. 

We are interested in the ramification filtrations of the decomposition 
groups acting on the completed local field at wild ramified points. 
We introduce the following notation:
We consider the set of jumps of the ramification filtration  $1=t_f < t_{f-1} < \cdots < t_1=n$, such that 
\begin{equation} \label{jumps-filt}
G_1=\ldots =G_{t_f}>G_{t_f+1}=\ldots=G_{t_{f-1}}>G_{t_{f-1}+1} \geq \ldots \geq G_{t_1}=G_{n}> \{1\},
\end{equation}
{\em i.e.}, $G_{t_i} > G_{t_i+1}$. For this sequence it is known that  $t_i\equiv t_j \mod p$
\cite[Prop. 10 p. 70]{SeL}.

\subsection{Lyndon-Hochschild-Serre Spectral Sequences}

In \cite{Hoch-Serre}, Hochschild and Serre considered the following 
problem: Given the short exact sequence of groups 
\begin{equation}\label{short-seq}
1\rightarrow H \rightarrow G \rightarrow G/H \rightarrow 1,
\end{equation}
and a $G$-module $A$, how are the cohomology groups 
\[H^i(G,A), H^i(H,A) \mbox{ and  } H^i(G/H,A^H)\]
related? They gave an answer to the above problem in 
terms of a spectral sequence. For small values of $i$ this spectral 
sequence gives us the low degree terms exact sequence:
\begin{equation}\label{low-LHS}
0\rightarrow H^1(G/H,A^H) \stackrel{{ \rm inf}}{\rightarrow}
H^1(G,A) \stackrel{{\rm res}}{\rightarrow} H^1(H,A)^{G/H} 
\stackrel{{\rm tg}}{\rightarrow} 
H^2(G/H,A^H) \stackrel{{\rm inf}}{\rightarrow} H^2(G,H),
\end{equation}
where $\mathrm{res},\mathrm{tg},\mathrm{inf}$ denote the 
restriction,transgression and inflation maps respectively.
%
 %
 % ------------------------------
 % Action of G/H on H^1(H,M)
 %--------------------------------
 %
 %
 \begin{lemma} \label{action-on-cocycles}
 Let $H$ be a normal subgroup of $G$, and let $A$ be a $G$-module. The group $G/H$ acts on the cohomology group $H^1(H,A)$ in terms 
 of the conjugation action given explicitly on the level of $1$-cocycles as follows:
 Let $\bar{\sigma}=\sigma H\in G/H$. The cocycle 
 \[
 \begin{array}{c}
 d:H \rightarrow A\\
 x \mapsto d(x)
 \end{array}
 \] 
 is sent by the conjugation action to the cocycle 
 \[
 \begin{array}{c}
 d^{\bar{\sigma}}:H \rightarrow A\\
 x \mapsto \sigma d(\sigma^{-1}  x \sigma) 
 \end{array},
 \] 
where $\sigma \in G$ is a representative of $\bar{\sigma}$.
 \end{lemma}
 \begin{proof}
 This  explicit description of the conjugation action on the level of cocycles is given in \cite[prop. 2-5-1, p.79]{Weiss}  The action is well defined by \cite[cor. 2-3-2]{Weiss}.
 \end{proof}

Our strategy is to use equation (\ref{low-LHS}) in order to reduce the problem 
of computation of $H^1(G,\TO)$ to an
 easier computation involving only 
elementary  abelian groups. 

\begin{lemma} \label{tame-computing}
	Let $A$ be a $k$-module, where $k$ is a field of characteristic $p$. 
	For the cohomology groups we have
	$H^1(G_0,A)=H^1(G_1,A)^{G_0/G_1}$.
\end{lemma}
\begin{proof}
Consider the short exact sequence 
	\[
	0 \rightarrow G_1 \rightarrow G_0 \rightarrow G_0/G_1 \rightarrow 0.
	\]
Equation (\ref{low-LHS}) implies the sequence
\[
0 \rightarrow H^1(G_0/G_1,A^{G_1}) \rightarrow H^1 (G_0,A) 
\rightarrow H^1(G_1,A)^{G_0/G_1} \rightarrow H^2(G_0/G_1,A^{G_1}).
\]
But the order of $G_0/G_1$ is not divisible by $p$, and is an invertible 
element in the $k$-module $A$. Thus the groups $
H^1(G_0/G_1,A^{G_1})$ and $H^2(G_0/G_1,A^{G_1})$ vanish 
 and the desired result follows \cite[Cor. 6.59]{Weibel}. 
\end{proof}

\begin{lemma} \label{trivial-action}
If $G=G_i,H=G_{i+1}$ are groups in the ramification filtration of 
the decomposition group at some wild ramified point, and $i\geq 1$ then 
the conjugation action of $G$ on $H$ is trivial.
\end{lemma}
\begin{proof}
     Let $L/K$ denote a wild ramified extension of local fields with Galois group $G$, let $\O_L$ denote the 
ring of integers of $L$ and let $m_L$ be the maximal ideal of $\O_L$. Moreover we will denote by $L^*$ the 
group of units of the field $L$.
	We can define \cite[Prop. 7 p.67,Prop 9 p. 69]{SeL} injections 
	\[
	\theta_0:\frac{G_0}{G_1} \rightarrow L^*
	\mbox{ and }
	\theta_i:
	\frac{G_i}{G_{i+1}} \rightarrow \frac{m_L^i}{m_L^{i+1}}, 
	\]
	with the property 
	\[
	\forall \sigma \in G_0 \mbox{ and } \forall \tau\in G_i/G_{i+1}:
	\theta_i(\sigma \tau \sigma ^{-1} ) =\theta_0 (\sigma)^i \theta_i(\tau)
	\]
	If $\sigma\in G_{t_j} \subset G_1$ then $\theta_0(\sigma)=1$ and 
	since $\theta_i$ is an injection, the above equation 
	implies that $\sigma \tau \sigma^{-1}=\tau$. 
	Therefore, the conjugation action of an element $\tau\in G_i/G_{i+1}$
	on $G_j$ is trivial, and the result follows. 
\end{proof}

\subsection{Description of the transgression map}
%
%  Alternative description of transgression
%
In this section we will try to determine the kernel of the transgression 
map.
The definition  of the transgression map given in (\ref{low-LHS}) is not suitable 
for computations. We will give an alternative description  following 
 \cite{NeuCoho}. 
 
Let $A$ be a $k$-algebra that is acted on by $G$ so that the $G$ action is compatible with the 
 operations on $A$.
 Let $\bar{A}$ be the set $\mathrm{Map}(G,A)$ of set-theoretic maps of the 
 finite group $G$ to the $G$-module $A$. The set $\bar{A}$ can be seen as a $G$-module by defining the action
 $f^g(\tau)=gf(g^{-1}\tau)$ for all $g,\tau \in G$. We observe that $\bar{A}$ is projective. The submodule $A$ can be seen as the 
 subset of constant functions.
Notice that the induced action of $G$ on the submodule $A$ seen as the submodule of  constant functions of $\bar{A}$ coinsides with the initial 
 action of $G$ on $A$.
 We consider the short exact sequence of $G$-modules:
 \begin{equation} \label{map-ses-ind}
 0 \rightarrow A \rightarrow \bar{A} \rightarrow A_1 \rightarrow 0.
 \end{equation}

 Let $H\lhd G$.  By applying the functor of $H$-invariants to 
the short exact sequence  (\ref{map-ses-ind})
 we obtain the long exact sequence 
 \begin{equation} \label{defpsi}
 0 \rightarrow A^H \rightarrow \bar{A}^H \rightarrow A_1^H
 \stackrel{\psi}{\rightarrow} H^1(H,A) \rightarrow 
 H^1(H,\bar{A})=0,
 \end{equation}
 where the last cohomology group is zero since $\bar{A}$ is projective. 
 
 We split the above four term sequence, by defining $B=\ker{\psi}$ to two 
 short exact sequences, namely:
 \[
 0 \rightarrow A^H \rightarrow \bar{A}^H \rightarrow B \rightarrow 0,
 \]
 \begin{equation} \label{defpsi2}
 0 \rightarrow B \rightarrow A_1^H \stackrel{\psi}{\longrightarrow} H^1(H,A) \rightarrow 0.
 \end{equation}
 Now we apply the $G/H$-invariant functor to the above two short exact 
 sequences in order to obtain:
 \[
 H^i(G/H,B)=H^{i+1}(G/H,A^H),
 \]
 and
 \begin{equation} \label{long-exact-ker-tg}
0 \rightarrow B^{G/H} \rightarrow A_1^G \rightarrow H^1(H,A)^{G/H} 
 \stackrel{\delta}{\rightarrow} H^1(G/H,B) \stackrel{\phi}{\rightarrow}  H^1(G/H, A_1^H) \cdots
\end{equation}
  It can be proved \cite[Exer. 3 p.71]{NeuCoho} that the composition
 \[
 H^1(H,A)^{G/H} \stackrel{\delta}{\rightarrow} H^1(G/H,B) \stackrel{\cong}{\rightarrow}
 H^2(G/H,A^H),
 \]
 is the transgression map. 
 \begin{lemma} \label{zero-trans}
Assume that $G$ is an abelian group. If the quotient $G/H$ is a cyclic group isomorphic to 
 $\mathbb{Z}/p\Z$ and the group $G$ can be written as a direct sum $G=G/H\times H$ then the transgression 
 map is identically zero. 
 \end{lemma}
 \begin{proof}
Notice that if $A^H=A$ then this lemma can be proved by the explicit form of the transgression 
map as a cup product \cite[Exer. 2 p.71]{NeuCoho}, \cite{Hoch-Serre}.

The study of the kernel of the transgression is reduced to the study of the kernel of $\delta$ 
 in (\ref{long-exact-ker-tg}). We will prove that the map $\phi$ in (\ref{long-exact-ker-tg}) is 
 1-1, and then the desired result will follow by exactness.
 
 Let $\sigma$ be a generator of the cyclic group $G/H=\mathbb{Z}/p\Z$.
 We denote by $N_{G/H}$ the norm map $A \rightarrow A$, sending 
\[A \ni a \mapsto \sum_{g\in G/H} ga= \sum_{\nu=0}^{p-1} \sigma^\nu a.\]
 By $I_{G/H}A$ we denote the submodule $(\sigma-1)A$ and by $_{N_{G/H}}A=\{a \in A: N_{G/H}a=0\}$.
 Since $G/H$ is a cyclic group we know that \cite[VIII 4]{SeL}, \cite[Th. 6.2.2]{Weibel}: 
 \[
 H^1(G/H,B)=\frac{_{N_{G/H}}B}{I_{G/H}B} \mbox{ and } 
 H^1(G/H,A_1^H)= \frac{_{N_{G/H}} A_1^H}{I_{G/H} A_1^H}.
 \]
 Thus, the map $\phi$ is given by 
 \[
 \frac{_{ N_{G/H}} B}{ I_{G/H} B} \rightarrow \frac{_{N_{G/H}} A_1^H}{ I_{G/H} A_1^H},
 \]
sending 
\[
b \mod I_{G/H} B \mapsto b \mod I_{G/H} A_1^H.
\]
The map $\phi$ is well defined since $I_{G/H} B \subset I_{G/H} A_1^H$.
 The kernel of $\phi$ is computed:
 \[
 \ker{\phi}=\frac{_{N_{G/H}} B \cap I_{G/H} A_1^H}{I_{G/H} B}.
 \]
The short exact sequence in (\ref{defpsi2}) 
 is a short exact sequence of $k[G/H]$-modules. This sequence seen as a short exact sequence 
of $k$-vector spaces is split, {\em i.e.} there is a section $s:H^1(H,A) \rightarrow A_1^H$ so that 
$\psi\circ s=\mathrm{Id}_{H^1(H,A)}$. This section map is only a $k$-linear map and not 
apriori compatible with the $G/H$-action.

Let us study the map $\psi$ more carefully. An element   $x\in A_1^H$ is a 
class $a \mod A$ where $a \in \bar{A}$, and since $x \in A_1^H$ we have that 
\[
a^h-a=ha-a =c[h] \in A.
\]
It is a standart argument that  $c[h]$ is an $1$-cocyle $c[h]:H\rightarrow A$ and the class
of this cocycle is defined to be $\psi(x)$.
Since the image of $c[h]$ seen as a cocycle $c[h]:H \rightarrow \bar{A}$ is trivial,
$c[h]$ is a coboundary {\em i.e.} we can select  $\bar{a}_c\in \bar{A}$ so 
that 
\begin{equation} \label{eqacdef}
c[h]= \bar{a}_c^h-\bar{a}_c.
\end{equation}
Obviously $\bar{a}_c \mod A$ is $H$-invariant and we define one section as
\[
s (c[h])= \bar{a}_c \mod A.
\]
We have assumed that the group $G$ can be written as 
 $G=H \times G/H$ therefore we can write the functions $\bar{a}_c$ as functions of two 
arguments
\[
 \bar{a}_c: \left\{
\begin{array}{ccc} 
 H \times G/H & \rightarrow & A\\
 (h,g) & \mapsto & \bar{a}_c(h,g)
\end{array}
\right.
,
\]
Notice that (\ref{eqacdef}) gives us that the for every $h,h_1\in H$ the 
 quantity $\bar{a}_c(h_1,g_1)^h-\bar{a}_c(h_1,g_1)$ does not depend on $g_1 \in G/H$.
% {\em i.e.} 
% \begin{equation} \label{equivcoc}
%  \bar{a}_c(h_
% \end{equation}
Now for any element $g\in G/H$ so that $g=\sigma H$ the action of $g$ on $c[h]$ is 
given by lemma \ref{action-on-cocycles}:
\begin{equation} \label{rr1}
c[h]^g=\sigma c[\sigma h \sigma^{-1}]=\sigma(  \bar{a}_c^{\sigma^{-1} h \sigma}-\bar{a}_c).
\end{equation}
The function $\bar{a}_c^{\sigma^{-1} h \sigma}$ is the function 
sending 
\begin{equation} \label{rr2}
H\times G/H \ni (h_1,g_1) \mapsto  \sigma^{-1} h \sigma \bar{a}_c  ( h^{-1}  h_1,  \sigma^{-1} \sigma g_1)=
\sigma^{-1} h \sigma \bar{a}_c  ( h^{-1}  h_1,  g_1)
\end{equation}
By combining (\ref{rr1}) and (\ref{rr2}) we obtain
\[
c[h]^g
=h \sigma \bar{a}_c -\sigma \bar{a}_c,
\]
since the function (notice the invariance on $g$ in the second argument)
\[
H \times G/H \ni  (h_1,g_1) \mapsto  h \sigma \bar{a}_c(h^{-1} h_1, g_1)- \sigma \bar{a}_c (h_1,g_1)=\]
\[=
  h \sigma \bar{a}_c(h^{-1} h_1, g^{-1} g_1)- \sigma \bar{a}_c (h_1, g^{-1} g_1)= \left( \bar{a}_c^{\sigma}\right)^h -
 \bar{a}_c^{\sigma}
\]
 The above proves that 
\[
s ( c[h]^\sigma ) =s (c[h])^\sigma,
\]
{\em i.e.}, the function $s$ is compatible with the $G/H$-action.
On the other hand,
every element $a\in A_1^H$ can be written as 
\[
a=b_a + s( \psi(a)),
\]
where $b_a:=a -s \psi(a) \in B$, since $\psi(b_a)=0$. 
The arbitary element in $I_{G/H} A_1^H$ is therefore written as
\begin{equation} \label{inter-df}
(\sigma-1)a = (\sigma-1)b_a + \sigma s \big(\psi(a) \big) - s \big(\psi(a) \big)=
\end{equation}
\[
=(\sigma-1)b_a+ s ( \sigma \cdot \psi(a) -\psi(a)).
\]
If $(\sigma-1)a \in  _{N_{G/H}} B \cap I_{G/H} A_1^H$ then since $\mathrm{Im}(s) \cap B =\{ 0\}$ we have that 
\[
s ( \sigma \cdot \psi(a) -\psi(a))=0 \Leftrightarrow (\sigma-1)a = (\sigma-1)b_a \in I_{G/H} B.
\]
Therefore, $\phi$ is an injection and the desired result follows.
\end{proof}
%
%
%

%\subsection{Reduction to the cyclic case}
%\input{redcyc.tex}
%

%
%
%
\subsection{The $G$-module structure of $\TO$.}
Our aim is to 
compute the first order infinitesimal deformations, 
{\em i.e.}, the tangent space  $D_\rho(k[\epsilon])$ to the infinitesimal 
deformation functor $D_\rho$ \cite[p.272]{MazDef}.  
This space can be  identified with $H^1(G,\TO)$. 
The conjugation action on  $\TO$ is defined as follows:
\begin{equation} \label{action-conj123}
\left(f(t)\frac{d}{dt}\right)^\sigma =f(t)^\sigma \sigma \frac{d}{dt}\sigma^{-1}
=
f(t)^\sigma \sigma\left( \frac{d\sigma^{-1} (t)}{dt} \right)\frac{d}{dt},
\end{equation}
where $\frac{d}{dt}\sigma^{-1}$ denotes the operator sending an element $f(t)$ to 
$\frac{d}{dt} f^{\sigma^{-1}}(t)$, {\em i.e.} we first compute the action of $\sigma^{-1}$ on $f$ and then 
we take the derivative with respect to $t$.
We will approach the cohomology group $H^1(G,\TO)$ using 
the filtration sequence given in (\ref{ram-filt}) and 
the low degree terms of the Lyndon-Hochschild-Serre spectral 
sequence.

The study of the cohomology group $H^1(G,\TO)$ can be reduced to 
the study of the cohomology groups $H^1(V,\TO)$, where $V$ is an 
elementary Abelian group. These groups can be written as a sequence of 
Artin-Schreier extensions that 
 have  the advantage that the extension 
and the corresponding actions have a relatively simple explicit 
form:

\begin{lemma}
Let  $L$ be a an elementary abelian $p$-extension  of the local field $K:=k((x))$, 
with Galois group $G=\oplus_{\nu=1}^s  \Z_p$, such that the maximal 
ideal of $k[[x]]$ is ramified completely and the ramification filtration has no intermediate  jumps {\em i.e.}  is given by 
\[
G=G_0=\cdots=G_n>  \{1\}=G_{n+1}.
\]
Then the extension $L$ is given  by $K(y_1,\ldots,y_s)$ where  
$1/y_i^p -1/y_i=f_i(x)$, where $f_i\in k((x))$ with a pole at the maximal ideal of order $n$. 
\end{lemma}
\begin{proof}
The desired result follows by the 
characterization of Abelian $p$-extensions in terms of Witt 
vectors, \cite[8.11]{Jac2}. Notice that the exponent of the group $G$ is $p$ and we 
have to consider  the image of $W_1(k((x))=k((x))$, where $W_\lambda (\cdot)$ denotes the Witt
ring of order $\lambda$ as is defined in \cite[8.26]{Jac2}.
\end{proof}

\begin{lemma} \label{artin-schreier-p}
Every $\Z/p\Z$-extension $L=K(y)$ of the local field $K:=k((x))$, 
with Galois group $G=  \Z/p\Z$, such that the maximal 
ideal of $k[[x]]$ is ramified completely, is given in terms of an  equation
$f(1/y)=1/x^n$, where $f(z)=z^p-z \in k[z]$.
The Galois group of the above extension can be identified with the $\mathbb{F}_p$-vector space 
$V$ of the roots of the polynomial $f$, and the correspondence is given by 
\begin{equation} \label{action_y}
\sigma_v: y\rightarrow \frac{y}{1+vy} \mbox{ for } v\in V.
\end{equation}
Moreover,   we can select a uniformization 
parameter of the local field $L$ such that the automorphism $\sigma_v$ 
acts on $t$ as follows:
\[
\sigma_v(t) =\frac{t}{(1+vt^n)^{1/n}}. 
\]
Finally, the ramification filtration is given by 
\[
G=G_0=\cdots=G_n>  \{1\}=G_{n+1},
\]
and $n\not \equiv 0 \mod p$.
\end{lemma}
\begin{proof}
By the characterization of Abelian extensions in terms of Witt 
vectors, \cite[8.11]{Jac2} we have that $f(1/y)=1/x^n$, where 
$f(z)=z^p-z\in k[z]$ (look also \cite[A.13]{StiBo}). 
Moreover the Galois group can be identified with the one dimensional $\mathbb{F}_p$-vector 
space $V$ of roots of $f$, sending $\sigma_v: y\rightarrow \frac{y}{1+vy}$.

The filtration of the ramification group $G$ is given by 
$G\cong G_0 =G_1=\cdots G_n,$ $G_i=\{1\}$ for $i\geq n+1$ 
\cite[prop. III.7.10 p.117]{StiBo}. 
By computation
\begin{equation} \label{form2.9}
x^n=((1/y)^p - 1/y )^{-1}=\frac{y^p}{1-y^{p-1} },
\end{equation}
hence $v_L(y)=n$, {\em i.e.}, $y=\epsilon t^n$, where $\epsilon$ is 
a unit in $\O_L$ and $t$ is the uniformization parameter in $\O_L$. 

Moreover, the polynomial $f$ can be selected so that 
 $p\nmid n$ \cite[III. 7.8.]{StiBo}. Since $k$ is an algebraically closed field, 
Hensel's lemma implies that every unit in $\O_L$ is an $n$-th 
power, therefore we might select the uniformization parameter $t$
such that $y=t^n$, and the desired result follows by (\ref{action_y}). 
\end{proof}

%
%-------------------------------------------------------------
%End of paste compute.tex
%
\begin{lemma} \label{upper-lower}
Let $H=\oplus_{\nu=1}^s \Z/p\Z$ be  an elementary 
Abelian group with ramification filtration 
\[
H=H_0=...=H_n> H_{n+1}=\{{\rm Id}\} \mbox{ and } H_\kappa=\{ \mathrm{Id} \}  \mbox{ for } \kappa \geq n+1.
\]
The upper ramification filtration in this case coincides with the lower ramification filtration. 
\end{lemma}
\begin{proof}
Let $m$ be a natural number. We define the function $\phi: [0,\infty] \rightarrow  \mathbb{Q}$ so that  for $m \leq u < m+1$ 
\[
\phi(u)=\frac{1}{|H_0|} \sum_{i=1}^m |H_i| + (u-m) \frac{ |H_{m+1}|}{|H_0|},
\]
and since $H_{n+1}=\{ \mathrm{Id} \}$ we compute 
\[
\phi(u)=
\left\{
\begin{array}{ll}
	u & \mbox{if } m+1 \leq n \\
	n+ \frac{u-n-1}{|H_0|} & \mbox{if } m+1 > n
\end{array}
	\right..
\]
The inverse function $\psi$ is computed by 
\[
\psi(u)=
\left\{
\begin{array}{ll}
	u & \mbox{if } u \leq n \\
	|H_0|u+ (-n |H_0|+n+1) & \mbox{if } u > n
\end{array}
\right..
\]
Therefore, by the definition of the upper ramification filtration we 
have $H^i=H_{\psi(i)}=H_i$ for $i \leq n$, while for $u>n$ we compute 
$\psi(u)=|H_0| u  - n |H_0|+n \geq n$, thus $H^u=H_{\psi(u)}=\{\mathrm{Id}\}$.
\end{proof}
\begin{lemma} \label{p-derivative}
Let $a\in \mathbb{Q}$. Then for every prime $p$ and every $\ell \in \mathbb{N}$ we have
\[
\lf \frac{ \lf \frac{a}{p^\ell} \rf }{p} \rf = \lf \frac{a}{p^{\ell+1}} \rf.
\]
\end{lemma}
\begin{proof}
Let us write the $p$-adic expansion of $a$:
\[
a= \sum_{\nu =\lambda} ^{-1} a_\nu p^\nu + \sum_{\nu=0}^\infty a_\nu p^\nu,
\]
where $\lambda \in \mathbb{Z}, \lambda<0$.
We compute 
\[
\lf \frac{a}{p^\ell} \rf =\sum_{\nu=\ell} ^\infty a_\nu p^\nu.
\]
and 
\[
\lf \frac{ \lf \frac{a}{p^\ell} \rf }{p} \rf =\sum_{\nu=\ell+1} ^\infty a_\nu p^\nu =\lf \frac{a}{p^{\ell+1}}\rf .
\]
\end{proof}
The arbitrary $\sigma_v \in \Gal(L/K)$ sends 
$t^n \mapsto \frac{t^n}{1+vt^n}$, so by computation 
\[
\frac{d\sigma_v(t)}{dt} =
\frac{1}
{
(1+vt^n)^{
           \frac{n+1}{n}
         }
}
\]
%{\n Remark:} 
%Let $x\in \mathbb{R}$. We will denote by  $\lf x \rf$ 
%the integer part of $x$, {\em i.e.}, the greatest  integer smaller than $x$ and by 
%$\lc x \rc$ the least integer greater than $x$. 
%In what follows we will made extensive  use of the operator $\lf \frac{\cdot}{p}\rf$. It is nice to point out that 
%the operator $\lf \frac{\cdot}{p}\rf$ behaves like a formal $p$-derivative on $\Z$. Indeed, let us 
%write an integer $n$ as $n=a_0 + a_1 p + \cdots + a_k p^k$, where $0 \leq a_i < p$. Then $\lf \frac{n}{p}\rf=a_1  + \cdots + a_k p^{k-1}$.
\begin{lemma}\label{action-desc}
We consider an  Artin-Schreier extension $L/k((x))$  and we keep the  notation from  lemma \ref{artin-schreier-p}.
Let  $\sigma_v\in \Gal(L/K)$. The corresponding action on the tangent space
$\TO$ is given by 
\[
\left(
f(t)\frac{d}{dt}
\right)^{\sigma_v} =f(t)^{\sigma_v} (1+vt^n)^{\frac{n+1}{n}}\frac{d}{dt}.
\]
\end{lemma}
\begin{proof}
We have that  
$\frac{d\sigma_v^{-1}(t)}{dt}=
\frac{d\sigma_{-v}(t)}{dt}=\frac{1}{(1-vt^n)^{\frac{n+1}{n}}}$ 
and by computation
\[
\sigma_v(\frac{d\sigma_{-v}(t)}{dt})=(1+vt^n)^{\frac{n+1}{n}}.
\] 
\end{proof}
Let $\O=\O_L$,
we will now compute the space of ``local modular forms'' 
\[
\TO^{G_{t_i}}=\{f(t)\in \O : f(t)^{\sigma_v}=f(t)(1+vt^n)^{-\frac{n+1}{n}}\},
\]for  $i\geq 1$.
First we do the computation for a cyclic $p$-group.
\begin{lemma} \label{artin-schreier-tangent}
Let $L/k((x))$ be an Artin-Schreier extension with Galois group  $H=\Z/p\Z$ and 
ramification filtration
\[
H_0=H_1=\ldots=H_n > \{\mathrm{Id} \}.
\]
Let $t$ be the uniformizer of $L$ and 
denote by $\TO$ the set of elements of the form $f(t)\frac{d}{dt}$, $f(t)\in k[[t]]$ equiped with the 
conjugation action defined in (\ref{action-conj123}) . 
The space $\TO^{G}$ is $G$-equivariantly isomorphic to 
the $\O_K$-module consisted of elements 
of the form 
\[
f(x)x^{n+1-\lf \frac{n+1}{p} \rf } \frac{d}{dx}, \;\;f(x)\in \O_K.
\]
\end{lemma}
\begin{proof}
Using the description of the action in lemma \ref{action-desc} we see
that
$\TO$ is isomorphic to the space of Laurent polynomials 
of the form $\{f(t)/t^{n+1}: f(t) \in \O\}$, and the isomorphism is 
compatible with the $G$-action. Indeed, we observe first that 
 $t^{n+1} \frac{d}{dt}$ is a $G$-invariant element in $\TO$. 
 Then, 
 for every 
$f(t)\frac{d}{dt}\in \TO$, the map sending 
\[
f(t)\frac{d}{dt}=\frac{f(t)}{t^{n+1}}t^{n+1}\frac{d}{dt} \mapsto 
\frac{f(t)}{t^{n+1}}, 
\]
is a $G$-equivariant isomorphism.

We have 
\[\{
f(t)/t^{n+1}, f(t) \in \O
\}^G=
\{
f(t)/t^{n+1}, f(t) \in \O
\}
\cap k((x)),
\]
so the $G$-invariant space 
consists of elements $g(x)$ in $K$ such that $g$ seen as an element 
in $L$ belongs to $\TO$, {\em i.e.}, $v_L(g) \geq -(n+1)$. 
Consider the set of functions $g(x)\in K$ such that 
$v_L(g)=p v_K(g)\geq -(n+1)$, {\em i.e.}, $v_K(g)\geq -\frac{n+1}{p}$.
Since $v_K(g)$ is an integer the last inequality is 
equivalent to $v_K(g)\geq \lc -(n+1)/p \rc = -\lf
(n+1)/p \rf$.

On the other hand, a simple computation with the defining equation 
of the Galois extension $L/K$ shows that 
\[
t^{n+1} \frac{d}{dt} =x^{n+1} \frac{d}{dx}, 
\]
and the desired result follows. 
\end{proof}
Similarly one can prove the more general:
\begin{lemma} \label{TO-invariants2-cyc}
We are using the notation of lemma \ref{artin-schreier-tangent}.
 Let $A$ be the fractional ideal $k[[t]]  t^a \frac{d}{dt}$ where $a$ is a fixed integer.
 The $G$-module $A$ is $G$-equivariantly isomorphic to $t^{a-(n+1)} k[[t]]$.
 Moreover, the space $A^G$  is the space of elements of the form 
 \[
 f(x) x^{n+1 - \lf \frac{n+1-a}{p} \rf} \frac{d}{dx}. 
 \]
\end{lemma}
Next we proceed to the more difficult case of ellementary abelian $p$-groups.
\begin{lemma} \label{TO-invariants1}
Let $G= \oplus_{i=1}^s \Z/p\Z$
 be the Galois group of the fully ramified elementary abelian extension $L/k((x))$
and assume that the ramification filtration is of the form 
\[
G=G_0=G_1=\cdots=G_n > \{1\}.
\]
Let $t$ denote the uniformizer  of $L$.
Denote by $\TO$ the set of elements of the form $f(t)\frac{d}{dt}$, $f(t)\in k[[t]]$ equiped with the 
conjugation action defined in (\ref{action-conj123}) . 
The space $\TO^{G}$ is $G$-equivariantly isomorphic to 
the $\O_K$-module consisted of elements 
of the form 
\[
f(x)x^{n+1-\lf \frac{n+1}{p^s} \rf } \frac{d}{dx}, \;\;f(x)\in \O_K,
\]
where $p^s=|G|$.
\end{lemma}
\begin{proof}
We will break the extension $L/k((x))$ to a sequence of extensions $L=L_0>L_1>\ldots L_s=k((x))$, 
such that $L_i/L_{i+1}$ is a cyclic $p$-extension. Denote by $\pi_i$ the uniformizer of $L_i$.
According to lemma \ref{upper-lower} the ramification extension $L_{i}/L_{i+1}$ is of conductor $n$, 
{\em i.e.} the conditions of  \ref{artin-schreier-tangent} are satisfied. 
We will prove the result inductively. For the extension $L/L_1$ the statement is true by lemma 
\ref{artin-schreier-tangent}. Assume that the lemma is true for $L/L_i$ so 
a $k[[\pi_i]]$ basis of $\TO^{\oplus_{\nu=1}^{i} \Z/p\Z}$ is 
given by the element $ \pi_i^{n+1-\lf \frac{n+1}{p^i} \rf } \frac{d}{d\pi_i}$.
Then lemma \ref{artin-schreier-tangent} implies that a $k[[\pi_{i+1}]]$ basis for 
\[
\TO^{\oplus_{\nu=1}^{i+1} \Z/p\Z}
=\left(
\TO^{
\oplus_{\nu=1}^{i} \Z/p\Z
}
\right)^{\Z/p\Z},
\]
is given by the element:
\[ \pi_{i+1}^{ n+1 -\lf \frac{n+1-\left(n+1-\lf \frac{n+1}{p^i} \rf \right)}{p} \rf } \frac{d}{d\pi_{i+1}}.
\]
The desired result follows by lemma \ref{p-derivative}.
\end{proof}
Similarly one can prove the more general:
\begin{lemma} \label{TO-invariants2}
We are using the notation of lemma \ref{TO-invariants1}.
 Let $A$ be the fractional ideal $k[[t]]  t^a \frac{d}{dt}$ where $a$ is a fixed integer.
 The $G$-module $A$ is $G$-equivariantly isomorphic to $t^{a-(n+1)} k[[t]]$.
 Moreover, the space $A^G$  is the space of elements of the form 
 \[
 f(x) x^{n+1 - \lf \frac{n+1-a}{p^s} \rf} \frac{d}{dx}. 
 \]
\end{lemma}

By induction, the above computation can be extended to the 
following:
\begin{pro}\label{comp-in-tan}
Let $L=k((t))$ be a local field acted on by a Galois $p$-group $G$ with ramification subgroups 
\[
G_1=\ldots=G_{t_f}> G_{t_f+1}=\ldots=G_{t_{f-1}}> G_{t_{f-1}+1}\geq \ldots \geq G_{t_1}=G_{n}> G_{t_0}=\{\mathrm{Id} \}.
\]
We consider the tower of local fields 
\[
L^{G_0}= L^{G_1}  \subseteq L^{G_i} \subseteq \ldots L^{\{\mathrm{Id} \}}=L.
\]
Let us denote by $\pi_i$ a local uniformizer for the field $L^{G_i}$, {\em i.e.} 
$L^{G_i}=k((\pi_i))$.   The extension 
$L^{G_{{t_i}+1}}/
L^{G_{t_i}}$ is Galois
with Galois group the elementary abelian group $H(i):=G_{t_i}/G_{{t_i}+1}$.
Moreover the ramification filtration of the group $H(i)$ is given by 
\[
H(i)_0=H(i)_1=\ldots=H(i)_{t_i} > H(i)_{t_i +1}=\{\mathrm{Id} \}
\]
and the conductor of the extension is $t_i$.
Let $\O$ be the ring of integers of $L$. 
The invariant space $\TO^{G_{t_i}}$ is the $\O^{G_{t_i}}$-module 
generated by:
\begin{equation}\label{compute-inv-tan}
\pi_i^{\mu_i} \frac{d}{d\pi_i}, 
\end{equation}
where $\mu_0=0$ and $\mu_i=t_i+1-\lf \frac{
-\mu_{i-1}+t_{i}+1
}{|G_{t_{i}}|/|G_{t_{i-1}}|} \rf$
\end{pro}
\begin{proof}
The first statements are clear from elementary Galois theory. What needs a proof is the formula for 
the dimensions $\mu_i$.
For $i=1$ we have that $G_{t_1}=G_n$ is an elementary abelian group and lemma 
\ref{TO-invariants1} applies, under the assumption $G_{t_0}=\{\mathrm{Id}\}.$ Therefore, 
\[
\TO^{G_{t_1}}=\pi_1 ^{
n+1 -\lf \frac{n+1}{|G_{t_1}|/|G_{t_0}|} \rf
}
\frac{d}{d\pi_1}.
\]
Assume that the formula is correct for $i$, {\em i.e.}, 
\[
\TO^{G_{t_i}}=\pi^{\mu_i} \frac{d}{d \pi_i}.
\]
Then lemma \ref{TO-invariants2} implies that 
\[
\TO^{G_{t_{i+1}}}=
\left(\TO^{G_{t_i}}\right)^{\frac{G_{t_{i+1}}}{{G_{t_i}}}}
=\pi_{i+1}^{\mu_{i+1}} \frac{d}{d\pi_{i+1}},
\]
where $\mu_{i+1}=n_{i+1} +1 - \lf \frac{n_{i+1}+1 -\mu_{i}  }{|G_{t_{i+1}}/{G_{t_i}}|} \rf$
and the inductive proof is complete. 
\end{proof}
%
%---------------------------------------------------------------------------
Let $k((t))/k((x))$ be a cyclic extension of local fields of order $p$,  
such that the maximal ideal $xk[[x]]$ is ramified completely in the 
above extension. 
For the ramification groups $G_i$ we have
\[
\Z/p\Z=G=G_0=\cdots =G_n > G_{n+1}=\{1\}. 
\]
Hence, the different exponent is computed $d=(n+1)(p-1)$. Let 
$E=t^a k[[t]]$ be a fractional ideal of $k((t))$. 
Let $N(E)$ denote the images of elements of $E$ under the norm map 
corresponding to the group $\Z/p\Z$. 
 It is known that  
$N(E)=x^{ \lf (d+a)/p \rf}k[[x]]$, and $E\cap k[[x]]=x^{\lc a/p \rc} k[[x]]$. 
The cohomology of cyclic groups is $2$-periodic and 
J. Bertin and A. M\'ezard in \cite[Prop. 4.1.1]{Be-Me}, proved that 
\begin{equation}\label{cyc-coh}
\dim_k H^1(G,E)=\dim_k H^2(G,E)=\frac{E\cap k[[x]]}{N(E)}=
\lf \frac{d+a}{p} \rf- \lc \frac{a}{p} \rc. 
\end{equation}
{\bf Remark:} The reader might notice that in \cite[Prop. 4.1.1]{Be-Me} instead of 
(\ref{cyc-coh}) the following formula is given:
\[
\dim_k H^1(G,k[[x]]\frac{d}{dx})=\lf \frac{2d}{p} \rf- \lc \frac{d}{p} \rc.
\]
But $k[[x]]\frac{d}{dx}\cong x^{-n-1}k[[x]]$, and $d=(n+1)(p-1)$, thus
\begin{eqnarray*}
\lf \frac{2d}{p} \rf- \lc \frac{d}{p} \rc= &\lf \frac{d + (n+1)p -n-1}{p} \rf - \lc \frac{ (n+1)p -n-1}{p} \rc= \\
= & \lf \frac{d  -n-1}{p} \rf - \lc \frac{  -n-1}{p} \rc, 
\end{eqnarray*}
and the two formulas coincide. 
\begin{cor} \label{dir-ab-groups}
Let $G$ be an abelian group that can be written as a direct product $G=H_1\times H_2$ of  groups 
$H_1,H_2$, and suppose that  $H_2=\Z/p\Z$. The following sequence is exact:
\[
0\rightarrow H^1(H_2,A^{H_1}) \rightarrow 
H^1(H_1\times H_2,A) \rightarrow H^1(H_1,A)^{H_2}
\rightarrow 0
\]
\end{cor}
\begin{proof}
 The group $H_2$ is cyclic of order $p$ so  the  transgression map is identically zero by 
  lemma \ref{zero-trans}  and the desired result follows. 
\end{proof}
{\bf Remark:}
It seems that the result of J. Bertin, A. M\'ezard, solves the problem of determining the 
dimension of the $k$-vector spaces  $H^1(\Z/p\Z,A)$ for fractional ideals of $k[[x]]$. But in what 
follows we have to compute the $G/H$-invariants of the above cohomology groups, therefore an 
explicit description of these groups and  of the $G/H$-action is needed.   
\section{Computing $H^1(\Z/p\Z,A)$. }
%
%
%
%Copy Paste to compute.tex 3Nov2004
%
We will need the following 
\begin{lemma} \label{binom-zero}
Let $a$ be a $p$-adic integer. The binomial coefficient $\binom{a}{i}$ is 
defined for  $a$ as usual:
\[
\binom{a}{i}=\frac{a(a-1)\cdot (a-i+1)}{i!}
\]
and it is also a $p$-adic integer \cite[Lemma 4.5.11]{Gou:97}.
Moreover, the binomial series is defined
\begin{equation}\label{bin-exp}
(1+t)^a = \sum_{i=0}^\infty \binom{a}{i} t^i. 
\end{equation}
Let $i$ be an integer and let $\sum_{\mu=0}^\infty b_\mu p^\mu$ and 
$\sum_{\mu=0}^\infty a_\mu p^\mu$
be the $p$-adic 
expansions  of $i$ and $a$ respectively. The $p$-adic integer 
$\binom{a}{i}\not \equiv 0 \mod p$ if and only if every 
coefficient $a_i \geq  b_i$. 
\end{lemma}
\begin{proof}
The only think that needs a proof is the criterion of the 
vanishing of the binomial coefficient $\mod p$. If $a$ is 
a rational integer, then this is a known theorem due to Gau\ss
\cite[Prop. 15.21]{Eisenbud:95}. When $a$ is a $p$-adic integer we compare the  
 coefficients $\mod p$ of the expression 
\[
(1+t)^a=(1+t)^{\sum_{\mu=0}^\infty  a_\mu p^\mu} =\prod_{\mu=1}^\infty 
(1+t^{p^\mu})^{a_\mu}
\]
and of the binomial expansion in (\ref{bin-exp})
and the result follows.
\end{proof}
\begin{lemma}[Nakayama map] \label{exp-Nakayama}
Let $G=\Z/p\Z$ be a cyclic group of order $p$ and let $A=t^a k[[t]]$ be a fractional ideal of $k[[t]]$.  Let $x$ be a local 
uniformizer of the field $k((t))^{\Z/p\Z}$.  
Let $\alpha \in H^2(G,A)$, and let $u[\sigma,\tau]$ be any cocycle representing the class $\alpha$.
The map 
\begin{equation} \label{exp-Nakayama2}
\phi:
H^2(G,A) \rightarrow  \frac{x^{\lc a/p \rc}k[[x]]}{x^{\lf \frac{(n+1)(p-1)+a}{p} \rf}k[[x]]}
\end{equation}
sending
\[
\alpha \mapsto \sum_{\rho\in G}u[\rho,\tau]\;\;\;\; \tau \in G,
\]
is well defined and is an  isomorphism. 
\end{lemma}
\begin{proof} 
Let $A$ be a $G$-module. Let us denote by  $\hat{H}^0(G,A)$ the zero Tate-cohomology. 
We use  Remark 4-5-7 and theorem 4-5-10 in the book of Weiss \cite{Weiss} in order to prove that the map 
$H^2(G,A) \ni \alpha \mapsto \sum_{\rho\in G}u[\rho,\tau] \in \hat{H}^0(G,A)$ is well defined and an isomorphism. 

Let $\sigma$ be a generator of the cyclic group $\Z/p\Z$.
We know that 
\[\hat{H}^0(\Z/p\Z,t^a k[[t]])=
\frac{\mathrm{ker}(\delta)}{N_{\Z/p\Z}( t^a k[[t]])},
\]
where $\delta=\sigma-1$ and $N_{\Z/p\Z}=\sum_{i=0}^{p-1} \sigma^i$. 
We compute that 
\[\ker(\delta)=t^ak[[t]] \cap k((x))=x^{\lc a/p\rc}k[[x]]\] and 
\[N_{\Z/p\Z}( t^a k[[t]])=x^{\lf \frac{ a+(n+1)(p-1)}{p} \rf} k[[x]]\] and this completes 
the proof. 
\end{proof}
Let $A=t^a k[[t]]$ be a fractional ideal of $k[[t]]$.  We consider the fractional ideal $t^{a+n+1}k[[t]]$, and we form the 
short exact sequence:
\begin{equation} \label{frac:ses1}
0 \rightarrow  t^{a+n+1} k[[t]] \rightarrow t^{a} k[[t]] \rightarrow M \rightarrow 0,
\end{equation} 
where $M$ is an $n+1$-dimensional $k$-vector space with basis 
$\{\frac{1}{t^{-a}},\frac{1}{t^{-a-1}},\ldots,\frac{1}{t^{-a-n}} \}.$

Let $\sigma_v$ be the automorphism $\sigma_v(t)=t /(1+vt^n)^{1/n}$, where $v\in \mathbb{F}_p$.
The action of $\sigma_v$ on $1/t^\mu$ is given by 
\begin{equation} \label{M-action1}
\sigma_v: \frac{1}{t^\mu} \mapsto \frac{(1+vt^n)^{\mu/n}}{t^\mu}=\frac{1}{t^\mu} 
\left(
\sum_{\nu=0}^\infty \binom{\mu/n}{\nu} v^\nu t^{\nu n}
\right).
\end{equation}
The action of $\Z/p\Z$ on the basis elements of $M$ is given by 
\begin{equation} \label{MZaction}
\sigma_v(1/t^\mu)=
\left\{
\begin{array}{ll}
1/t^\mu & \mbox{if } -a<\mu \\
1/t^{-a} +\frac{-a}{n} v 1/t^{-a-n}  & \mbox{if } \mu=-a 
\end{array}
\right. 
\end{equation}
We consider the long exact sequence we obtain by applying the  $G$-invariants functor on  (\ref{frac:ses1}):
\begin{equation} \label{long:ses1}
0 \rightarrow t^{a+n+1} k[[t]]^G \rightarrow t^{a} k[[t]]^G \rightarrow M^G  \stackrel{\delta_1}{\rightarrow}
\end{equation}
\[
\rightarrow H^1(G, t^{a+n+1} k[[t]]) \rightarrow H^1(G,t^{a} k[[t]]) \rightarrow H^1(G,M)\stackrel{\delta_2}{\rightarrow}
 H^2(G, t^{a+n+1} k[[t]]) \rightarrow \cdots
\]
\begin{lemma} \label{last-delta1}
Assume that the group $G=\Z/p\Z$ generated by $\sigma_v$.
The map $\delta_1$ in (\ref{long:ses1}) is onto. 
\end{lemma}
\begin{proof}
By (\ref{MZaction}) we have: 
\[
\dim_k M^{\Z/p\Z}=
\left\{
\begin{array}{ll}
n+1 & \mbox{if } p \mid a \\
n & \mbox{if } p \nmid a
\end{array}
\right..
\]
On the other hand, if $x$ is a local uniformizer of the field $k((t))^{\Z/p\Z}$,  then:
\[
\left(
\frac{1}{t^{-a-(n+1)}}k[[t]]
\right)^{\Z/p\Z}=
x^{\lc  \frac{a  +(n+1) }{p}\rc}k[[x]]=\frac{1}{x}^{\lf  \frac{-a-(n+1)}{p}\rf}k[[x]],
\]
and similarly 
\[
\left(
\frac{1}{t^{-a}}k[[t]]
\right)^{\Z/p\Z}=
\frac{1}{x}^{\lf \frac{-a}{p} \rf} k[[x]].
\]
The image of $\delta_1$ has dimension:
\[
 \dim_k M^{\Z/p\Z}-\lf \frac{-a}{p} \rf + \lf \frac{-a -(n+1)}{p} \rf.
\]
Moreover for the dimension of $H^1(\Z/p\Z, \frac{1}{t^{-a -(n+1)}} k[[t]])$ we compute:
\[
h:=\dim_k H^1(\Z/p\Z, \frac{1}{t^{-a -(n+1)}} k[[t]])=\]
\[
=\lf \frac{(n+1)(p-1)+a+(n+1)}{p} \rf - \lc \frac{a + (n+1)}{p}\rc=
\]
\[
=(n+1) - \lc \frac{-a}{p} \rc + \lf  \frac{-a -(n+1)}{p} \rf.
\]
We distinguish the following two cases:
\begin{itemize}
\item 
If $p \mid a$ then $\lf \frac{-a}{p} \rf= \lc \frac{-a}{p} \rc$, $\dim_k M^G=n+1$ and we observe that $\dim_k\mathrm{Im}(\delta_1)=h$.
\item 
If $p \nmid a$ then $\lc \frac{-a}{p} \rc = \lf \frac{-a}{p} \rf +1$, and $\dim_k M^G=n$ and in this case it also holds $\dim_k\mathrm{Im}(\delta_1)=h$.
\end{itemize} 
\end{proof}
\begin{pro} \label{h1M}
  The cohomology group  $H^1(\Z/p\Z,M)$ is isomorphic to 
  \[
  H^1(\Z/p\Z,M)\cong
  \left\{
  \begin{array}{ll}
 \bigoplus_{i=-a-n}^{-a} Hom(\Z/p\Z,k) & \mbox{if } p\mid a \\
 \bigoplus_{i=-a-n+1}^{-a} Hom(\Z/p\Z,k) & \mbox{if } p\nmid a
\end{array}
  \right..
  \] 
  \end{pro}
  \begin{proof}
Assume that the arbitrary automorphism $\sigma_v\in G=\Z/p\Z$ is given by $\sigma_v(t)=t/(1+vt^n)^{1/n}$ where 
$v\in \mathbb{F}_p$.

  Let us write a cocycle $d$ as 
  \[
  d\sigma_v=\sum_{i=-a-n}^{-a} \alpha_i(\sigma_v)\frac{1}{t^i}.
  \]
  By computation, 
  \[
  d(\sigma_v)^{ \sigma_v^\mu} = \sum_{i=-a-n}^{-a} \alpha_i(\sigma_v)\frac{1}{t^i}+
  \alpha_{-a}(\sigma_v)\frac{-a}{n}\mu v \frac{1}{t^{-a-n}}.
  \]

  Moreover,  the cocycle condition $d(\sigma_v+\sigma_w)=d(\sigma_w)+d(\sigma_v)^{\sigma_w}$ for $d(\sigma_v)=\sum_{i=-a-n}^{-a} \alpha_i(\sigma_v) \frac{1}{t^i}$  gives:
  \[
  \sum_{i=-a-n}^{-a} \alpha_i(\sigma_v+\sigma_w) \frac{1}{t^i}=
  \left(
  \sum_{i=-a-n}^{-a} \alpha_i(\sigma_v) \frac{1}{t^i}
  \right)^{\sigma_w} + 
  \sum_{i=-a-n}^{-a} \alpha_i(\sigma_w) \frac{1}{t^i}=\]
  \[
  \sum_{i=-a-n}^{-a} \alpha_i(\sigma_v) \frac{1}{t^i}
   + \alpha_{-a}(\sigma_v)\frac{-a}{n}w \frac{1}{t^{-a-n}}+
  \sum_{i=-a-n}^{-a} \alpha_i(\sigma_w) \frac{1}{t^i}
  \]
  By comparing coefficients we obtain:
  \[
  \alpha_i(\sigma_w+\sigma_v)=\alpha_i(\sigma_w)+\alpha_i(\sigma_v) \;\;\mbox{for } i \neq -a-n.
  \]
  and 
  \[
  \alpha_{-a-n}(\sigma_w+\sigma_v)=\alpha_{-a-n}(\sigma_w)+\alpha_{-a-n}(\sigma_v)+ \alpha_{-a}(\sigma_v)\frac{-a}{n} w. 
  \]
  The last equation allows us to compute the value of $\alpha_{-a-n}$ on any power $\sigma_v^\nu$ of the generator $\sigma_v$ of 
  $\Z/p\Z$. Indeed, we have:
  \[
  \alpha_{-a-n}(\sigma_v^\nu)=\nu \alpha_{-a-n}(\sigma_v) + (\nu-1)\alpha_{-a}(\sigma_v)\frac{-a}{n} v.
  \]
  This proves that the function  $\alpha_{-a-n}$ depends only on the selection of $\alpha_{-a-n}(\sigma_v)\in k $.
  
  We will now compute the coboundaries. Let $b=\sum_{i=-a-n}^{-a} b_i \frac{1}{t^i}$, $b_i \in k$ be an element in $M$. 
 By computation, 
 \[
 b^{\sigma_v}-b=b_{-a} \frac{-a}{n} v\frac{1}{t^{-a-n}}.
 \]
 We distinguish the following cases:
 \begin{itemize}
 \item 
 If $p\mid a$ then the $\Z/p\Z$ action on $M$ is trivial, so 
 \[H^1(\Z/p\Z,M)=Hom(\Z/p\Z,M) =\bigoplus_{i=-a-n}^{-a}  Hom(\Z/p\Z,k).\]
 The dimension of $H^1(\Z/p\Z,M)$ in this case is $n+1$.
 \item 
 If $p\nmid a$, then the coboundary kills the contribution of the cocycle on the $\frac{1}{t^{-a-n}}$ basis element and the
 cohomology group is 
  \[H^1(\Z/p\Z,M)=Hom(\Z/p\Z,M) =\bigoplus_{i=-a-n+1}^{-a}  Hom(\Z/p\Z,k).\]
 \end{itemize}
 \end{proof}

\begin{lemma} \label{last-delta2}
Assume that $p\geq 3$.
Let $e=1$ if $p\nmid a$ and $e=0$ if $p\mid a$.
If $n\geq 2$ then 
an element 
\[\sum_{i=-a-n+e}^{-a} a_i(\cdot) \frac{1}{t^i}\in H^1(\Z/p\Z, M)\]
 is in the kernel of 
$\delta_2$ if and only if $a_i(\cdot) \binom{i/n}{p-1}=0$ for all $i$.
If $n=1$ then an element 
\[\sum_{i=-a-n+e}^{-a} a_i(\cdot) \frac{1}{t^i}\in H^1(\Z/p\Z, M)\]
is in the kernel of $\delta_2$ if and only if 
$a_i(\cdot) \binom{i/n}{p-1}=0$ for all $-a-n+e\leq i \leq -a$ and 
$a_i(\cdot) \binom{i/n}{2p-2}=0$ for all $-a-n+e\leq i \leq -a$ such that 
$2(p-1)n-i < (n+1)p+ p \lf \frac{a}{p} \rf.$
\end{lemma}
\begin{proof}
A derivation $a_i(\sigma_v) \frac{1}{t^i}$, $ -a-n+e\leq i \leq -a$ representing a cohomology class in $H^1(\Z/p\Z,M)$ 
is mapped to 
\begin{eqnarray} \label{last-delta22} 
\delta_2 \big( a_i(\cdot) \frac{1}{t^i}\big)[\sigma_w,\sigma_v]  = &
a_i(\sigma_v) \frac{1}{t^i}^{\sigma_w} - 
a_i(\sigma_v+\sigma_w) \frac{1}{t^i} +
a_i(\sigma_w) \frac{1}{t^i}= \\
  =  & \frac{a_i(\sigma_v)}{t^i} 
  \left(\sum_{\nu=1}^\infty\binom{i/n}{\nu} w^\nu t^{n\nu} \right)  
              \nonumber
\end{eqnarray}
We consider now the map $\phi$ defined in (\ref{exp-Nakayama2}) in the proof of lemma \ref{exp-Nakayama}. 
The map $\delta_2:H^1(G,M) \rightarrow H^2(G,t^{a+n+1} k[[t]])$ is composed with $\phi$ and the 
the image of $\phi \circ \delta_2$ in ${x^{\lc \frac{a+n+1}{p} \rc}k[[x]]}/{x^{\lf \frac{(n+1)p+a}{p} \rf}k[[x]]}$ is given by 
\[
\phi\circ \delta_2 \big( a_i(\cdot) \frac{1}{t^i}\big)=
\sum_{w\in \Z/p\Z}
\frac{a_i(\sigma_v)}{t^i} 
  \left(\sum_{\nu=1}^\infty\binom{i/n}{\nu} w^\nu t^{n\nu} \right).
\]
On the other hand  recall that 
\[
\sum_{w\in \Z/p\Z} w^\nu=
\left\{
\begin{array}{rl}
0 & \mbox{if }  p-1 \nmid \nu\\
-1 & \mbox{if } p-1 \mid \nu
\end{array}
\right..
\]
and every homomorphism $a_i : (\Z/p\Z,\cdot) \rightarrow (k,+)$ is given by 
$a_i(\sigma_w)=\lambda_i w$, where $\lambda_i\in k$. 
Therefore, 
\[
\phi\circ \delta_2 \big( a_i(\cdot) \frac{1}{t^i}\big)
=\sum_{\nu=1}^\infty
\binom{i/n}{\nu} 
\left(
\sum_{w\in \Z/p\Z} w^\nu \right) a_i(\sigma_v)
t^{n\nu-i}=
\]
\begin{equation} \label{comp-last32}
\sum_{\nu=1,p-1 \mid \nu }^\infty
\binom{i/n}{\nu} 
(-1) a_i(\sigma_v)
t^{n\nu-i}.
\end{equation}
Observe that $p-1 \mid \nu$ is equivalent to $\nu=\mu p -\mu $, and since $\nu \geq 1$, we have $\mu \geq 1$. 
Thus, equation (\ref{comp-last32}) becomes
\[ 
  \sum_{\mu=1}^\infty\binom{i/n}{\mu p -\mu  } (-1) \lambda_i t^{(\mu p-\mu )n-i}= 
 \]
 \[
 =\binom{i/n}{ p -1 } (-1) \lambda_i t^{(p-1)n-i} +  \binom{i/n}{2 p -2 } (-1) \lambda_i t^{(2p-2)n-i} +\mbox{higher order terms}
 \]
%
% dior8wsh 25/10/2005 h arxikh eksiswsh htan la8os
%
 {\bf Claim:}  If $n\geq 2$ and $p\geq 3$ then for all   $a\leq -i \leq a +n$ and for $\mu \geq 2$ 
\begin{equation} \label{8877}
\mu(p-1)n-i  \geq p \lf  \frac{(n+1)p +a}{p}\rf.
\end{equation}
If $n=1$ and $p\geq 3$ then (\ref{8877}) holds for $a\leq -i \leq a +n$ and for $\mu \geq 3$.
Moreover 
\[
(p-1)n-i < p \lf  \frac{(n+1)p +a}{p}\rf,
\]
for $a \leq -i \leq a+n$.

Indeed, the inequality 
\begin{equation} \label{8878}
\mu \geq \frac{n+1}{n}\frac{p}{p-1}, 
\end{equation}
holds for $p\geq 3, n \geq 2$ and $\mu \geq 2$ or for $p \geq 3,n=1, \mu\geq 3$.
Therefore, (\ref{8878}) implies that 
\[
n+1 + \lf \frac{a}{p}\rf -\frac{a}{p} \leq n+1 \leq \mu \frac{p-1}{p}n \Rightarrow 
\]
\[
(n+1)p + \lf \frac{a}{p} \rf p \leq \mu (p-1)n+a\leq \mu (p-1)n-i 
\]
and the first assertion is proved.
On the other hand
\[
\frac{a}{p} < 1 + \lf \frac{a}{p} \rf \Rightarrow  \frac{a}{p} + n < n+1 + \lf \frac{a}{p} \rf  \Rightarrow
\]
\[
(p-1)n-i < a+ pn< p(n+1)+  p\lf \frac{a}{p} \rf,
\]
and the second assertion is proved.

Since for elements $g\in k[[x]] \subset k[[t]]$ we have $p v_x(g) = v_t(g)$ we observe that
all elements in $k[[t]]$ that have valuation greater or equal to $(n+1)p + \lf \frac{a}{p} \rf$ are zero in 
the lift of the ideal $x^{(n+1+ \lf \frac{a}{p} \rf} k[[x]]$ on $k[[t]]$. 
Therefore the claim gives us that that for $p\geq 3, n\geq 2$,
 \[\phi \circ \delta_2( a_i(\cdot) \frac{1}{t^i})= \binom{i/n}{ p -1 } (-1) \lambda_i t^{(p-1)n-i} \]
 so 
 $\sum_{i=-a-n}^{-a}  a_i(\cdot) \frac{1}{t^i}$ is in the kernel of $\delta_2$ if and only if 
 \[\binom{i/n}{ p -1 } (-1) \lambda_i=0\;\; \mbox{for all } i.\]
The case $n=1$ follows by a similar argument.
  \end{proof}

 \begin{pro}  \label{final-cyc-comp}
 The cohomology group $H^1(\Z/p\Z,t^a k[[t]])$ is isomorphic to the $k$-vector space generated by  
 \[
 \left\{
 \frac{1}{t^i}, b\leq i \leq -a, \mbox{ such that } \binom{\frac{i}{n}}{p-1}=0 
 \right\},
 \]
 where $b=-a-n$ if $p \mid a$ and $b=-a-n+1$ if $p\nmid a$. 
 \end{pro}
 \begin{proof}
If $n\geq 2$ and $p\geq 3$ then the 
  result is immediate  by the exact sequence (\ref{long:ses1}), lemmata \ref{last-delta1} and \ref{last-delta2}  and by the computation of $H^1(\Z/p\Z,M)$ 
 given in 
 proposition \ref{h1M}. 
 
Assume that $n=1$, and let $a=a_0 +a_1p + a_2p^2 +\cdots $ be the $p$-adic expansion of $a$.
Then the inequality 
\begin{equation} \label{123455}
(n+1)p + p\lf \frac{a}{p} \rf \leq 2 (p-1)n+a
\end{equation}
holds if $a_0\neq 0,1$. Indeed, in this case we have that $\frac{2}{p} \leq  \frac{a}{p}-\lf \frac{a}{p} \rf  <1$ and 
(\ref{123455}) holds.
Therefore, for the case $p\mid a$ and $a=1+p b$, $b\in \Z$ we have to check the binomial coefficients $\binom{\frac{i}{n}}{2p-2}$ as well.
We will prove that in these cases if $\binom{\frac{i}{n}}{p-1}=0$ then $\binom{\frac{i}{n}}{2p-2}=0$ and the proof will be complete.

Assume, first that $p\mid a$ and $n=1$. Then,$-a-1 \leq i \leq -a$, {\em i.e.} $i=-a-1$ or $i=-a$.
 We compute that $\binom{-a}{p-1}=0$ since there is 
no constant term in the $p$-adic expansion of $-a$. Moreover the $p$-adic expansion of $2p-2$ is computed 
$2p-2=p-2 + p$, and since $p\neq 2$ we have that $\binom{-a}{2p-2}=0$ as well.
For $i=-a-1$ we have that  $i=p-1 + p b$ for some $b\in \Z$ therefore by comparing the $p$-adic expanesions
of $-a-1$,$p-1$ we obtain that $\binom{-a-1}{p-1} \neq 0$, and this value of $i$ does not contribute to the cohomology.

Assume now that $a=1  + pb$, $b\in \Z$. We have that 
$i=-a$, and $-a=p-1 + p(b+1)$. Therefore by comparing the $p$-adic expansions of $-a,p-1$ we obtain that 
$\binom{-a}{p-1}\neq 0$ and this value of $i$ does not contribute to the cohomology.

\end{proof}

\begin{pro} \label{triv11}
Let $A=t^a k[[t]]$ be a fractional ideal of the local field $k((t))$.
Assume that $H=\oplus_{\nu=1}^s \Z/p\Z$ is an elementary 
Abelian group with ramification filtration 
\[
H=H_0=...=H_n> H_{n+1}=\{{\rm Id}\}.
\]
Let $\pi_i$ be the local uniformizer of the local field 
$k((t))^{\oplus_{\nu=1}^{i-1} \Z/p\Z}$, and 
$a_i=\lc \frac{a_{i-1}}{p} \rc$, $a_1=a$. 
The cohomology group $H^1(H,A)$ is generated as a $k$-vector space by the 
following basis elements:
\[
\left\{
\bigoplus_{\lambda=1}^s
\frac{1}{\pi_i^{i_\lambda}},\;\; 
\begin{array}{c} \lambda=1,\ldots,s \\
b_i \leq i_\lambda \leq  -a_i 
\end{array}
\mbox{ such that } \binom{i_\lambda/n}{p-1}=0,
\right\}
\] 
where $b_i=-a_i-n$ if $p \mid a_i$ and $b_i=-a_i-n+1$ if $p\nmid a_i$.
Moreover, let $H(i):=H/\oplus_{\nu=1}^{i-1} \Z/p\Z$.  
The groups $H^1(\oplus_{\nu=1}^{i-1} \Z/p\Z,t^a k[[t]])$ 
are trivial $H(i)$-modules
with respect to the conjugation action.
\end{pro}
\begin{proof}
For $A=t^a k[[t]]$, we compute the invariants 
$t^a k[[t]] \cap k((t))^{\Z/p\Z}=x^{\lc \frac{a}{p}\rc} k[[x]]$, 
where $x$ is a local uniformizer for the ring of integers 
of $k((t))^{\Z/p\Z}$. 

The modules $A^{\oplus_{\nu=1}^{i-1} \Z/p\Z}$ can be computed 
recursively:
\[
A^{\oplus_{\nu=1}^{i-1} \Z/p\Z}=\pi_i^{a_i} k[[\pi_i]],
\]
where $\pi_i$ is a uniformizer for the local field 
$k((t))^{\oplus_{\nu=1}^{i-1} \Z/p\Z}$ and $a_i=\lc \frac{a_{i-1}}{p} \rc$,
$a_1=a$.

In order to compute the ramification filtration of quotient groups we 
have to employ the upper ramification filtration for the ramification group 
\cite[IV 3, p. 73-74]{SeL}.
But according to lemma \ref{upper-lower} the upper ramification filtration coincides with 
the lower ramification filtration therefore
the ramification filtration for the groups 
$H(i)$ is
\[
H(i)_0=...=H(i)_n > \{{\rm Id}\}.
\]
For the group  $\frac{\Z}{p\Z} \oplus \frac{\Z}{p\Z}$  corollary \ref{dir-ab-groups} implies that 
\[
H^1({\Z}/{p\Z} \oplus {\Z}/{p\Z}, t^a k[[t]])= 
H^1\left(\frac{\frac{\Z}{p\Z} \oplus \frac{\Z}{p\Z}}{\Z/p\Z}, t^a k[[t]] ^{\Z/p\Z}\right) \oplus H^1(\Z/p\Z,t^a k[[t]])^{\frac{H}{\Z/p\Z}}.
\]
By lemma \ref{upper-lower} and by the compatibility  \cite[IV 3, p. 73-74]{SeL} of the upper ramification filtration with quotients, we obtain that the 
 quotient $\frac{H}{\Z/p\Z}$ has also conductor $n$. 
 By lemmata \ref{action-on-cocycles}, \ref{trivial-action} and by the explicit description of the group $H^1(\Z/p\Z,t^a k[[t]])$ of 
 proposition \ref{final-cyc-comp} and by the fact $\frac{H}{\Z/p\Z}$ is of conductor $n$, the action of $\frac{H}{\Z/p\Z}$ on $H^1(\Z/p\Z,t^a k[[t]])$ is trivial.
Thus 
\[
H^1(\Z/p\Z \oplus \Z/p\Z, t^a k[[t]])= H^1\left(\frac{H}{\Z/p\Z}, t^a k[[t]] ^{\Z/p\Z}\right) \oplus H^1(\Z/p\Z,t^a k[[t]]).
\]
Moreover the cohomology group $H^1(\Z/p\Z \oplus \Z/p\Z, t^a k[[t]])$ is generated over $k$ by 
$
\left\langle 
\frac{1}{\pi_1^i}  \bigoplus  \frac{1}{\pi_2^j}, 
\right\rangle
$
where  $b_1\leq i \leq -a$, $b_2 \leq j \leq  -\lc \frac{a}{p}\rc $  
and  $\binom{i/n}{p-1}=\binom{j/n}{p-1}=0$.

The desired result follows by induction.
%%%%%%%%--------->
\end{proof}
\begin{pro} \label{comp-el-ab}
Let $A=t^a k[[t]]$ be a fractional ideal of the local field $k((t))$.
Assume that $H=\oplus_{\nu=1}^s \Z/p\Z$ is an elementary 
Abelian group with ramification filtration 
\[
H=H_0=...=H_n> H_{n+1}=\{{\rm Id}\}.
\]
The dimension of $H^1(H,A)$ can be computed as
\begin{equation}\label{eq-el-ab}
\dim_k H^1(H,A)=\sum_{i=1}^s \left( \lf \frac{(n+1)(p-1)+a_i}{p} \rf
-\lc \frac{a_i}{p} \rc \right),  
\end{equation}
where $a_i$ are defined recursively by $a_1=a$ and $a_i=\lc \frac{a_{i-1}}{p} \rc$.
In particular if $A=k[[t]]$, then 
\begin{equation}\label{cyc-coh2}
\dim_kH^1(H,k[[t]])=s \lf  \frac{(n+1)(p-1)}{p} \rf. 
\end{equation}
\end{pro}
\begin{proof}
By induction on the number of direct sumands,corollary \ref{dir-ab-groups} 
and proposition \ref{triv11} we can prove the following formula:
\begin{equation}\label{gen-fo}
H^1(H,A)=\oplus_{i=1}^s H^1(\Z/p\Z,A^{\oplus_{\nu=1}^{i-1} \Z/p\Z}). 
\end{equation}

In order to compute the dimensions of the direct sumands 
$H^1(\Z/p\Z,A^{\oplus_{\nu=1}^{i-1} \Z/p\Z})$, for various $i$ 
we have to compute the ramification filtration for the groups defined
as  $H(i)=\frac{H}{\oplus_{\nu=1}^{i-1} \Z/p\Z}$, since 
$\oplus_{\nu=1}^{i-1} \Z/p\Z=H/H(i)$.
But the upper ramification filtration coincides with the lower 
ramification filtration \ref{upper-lower}.
Thus, the dimension of $H^1(H,A)$ can be computed as
\begin{equation}
\dim_k H^1(H,A)=\sum_{i=1}^s \lf \frac{(n+1)(p-1)+a_i}{p} \rf
-\lc \frac{a_i}{p} \rc.  
\end{equation}
In particular if $A=k[[t]]$, then 
\begin{equation}
\dim_kH^1(H,k[[t]])=s \lf  \frac{(n+1)(p-1)}{p} \rf. 
\end{equation}
\end{proof}

%
%    Bound approach
%
For the dimension 
\begin{equation}\label{kappa-def}
\kappa_{t_i}=\dim_k\ker\left(\tg : H^1(G_{t_i+1}, 
\TO) \rightarrow H^2(G_{t_{i}}/G_{t_i+1}, \TO^{G_{t_i+1}} )\right) 
\end{equation}
of the kernel of the transgression map the following formula holds:
\begin{equation}\label{upper-bound1}
0 \leq \kappa_i \leq  \dim_k  H^1(G_{t_i+1}, \TO)^{G_{t_{i}}/G_{t_i+1}} \leq \dim_k  H^1(G_{t_i+1},\TO). 
\end{equation}
This allows us to compute:
\begin{pro}\label{loc-def-pro}
Let $G$ be the Galois group of the extensions of local fields
$L/K$, with ramification filtration $G_i$ and let $(t_\lambda)_{1\leq \lambda \leq f}$ be the jump sequence in (\ref{jumps-filt}).
For the dimension of $H^1(G_1,\TO)$ the following bound holds:
\begin{eqnarray}\label{loc-def-dim}
  H^1(G_1/G_{t_{f-1}},\TO^{G_{t_{f-1}}}) \leq &  \dim_k H^1(G_1,\TO) \leq  \\
  \leq & \sum_{i=1}^f \dim_k H^1(G_{t_i}/G_{t_{i-1}},\TO^{G_{t_{i-1}}})^{\frac{G_{t_f}}{G_{t_i}}}  \leq \nonumber \\
  \leq & \sum_{i=1}^f \dim_k H^1(G_{t_i}/G_{t_{i-1}},\TO^{G_{t_{i-1}}}), \nonumber
\end{eqnarray}
where $G_{n+1}=\{1\}$.
The left bound is best possible 
in the sense that there  are ramification 
filtrations such that  the first  inequality
becomes equality. 
\end{pro}
\begin{proof}
Using the low-term sequence in (\ref{low-LHS}) we obtain the following
inclusion for $i\geq 1$:
\begin{equation} \label{ind-33}
H^1(G_{t_i},\TO)= H^1(G_{t_{i}}/G_{t_{i-1}}, \TO^{G_{t_{i-1}}}) +
 \ker{tg} \subseteq 
  \end{equation}
 \[
 \subseteq H^1(G_{t_{i}}/G_{t_{i-1}}, \TO^{G_{t_{i-1}}})\oplus H^1(G_{t_{i-1}},\TO)^{G_{t_{i}}/G_{t_{i-1}}}.
\]
We start our computation from the end of the ramification groups:
\begin{equation} \label{ind-11}
H^1(G_{t_2},\TO) \subseteq H^1(G_{t_2}/{G_{t_1}},\TO^{G_{t_1}}) \oplus H^1(G_{t_1},\TO)^{G_{t_2}/{G_{t_1}}}.
\end{equation}
Observe here that $\TO$ is not $G_{t_1}$-invariant so there in no apriori well defined action of $G_{t_2}/{G_{t_1}}$ on $\TO$. 
But since the group $G_{t_1}$ is of conductor $n$  using the explicit form of $H^1(G_{t_1},\TO)$ we see that 
that $H^1(G_{t_1},\TO)$ is a trivial $G_1$-module. Of course this is also clear from the general properties of the 
conjugation action \cite[cor. 2-3-2]{Weiss}.
We go on to the next step:
\begin{equation} \label{ind-22}
H^1(G_{t_3},\TO)\subseteq H^1(G_{t_3}/G_{t_2},\TO^{G_{t_2}})\oplus H^1(G_{t_2},\TO)^{G_{t_3}/G_{t_2}}.
\end{equation}
The combination of (\ref{ind-11}) and (\ref{ind-22}) gives us 
\[
H^1(G_{t_3},\TO) \subseteq H^1\left(\frac{G_{t_3}}{G_{t_2}},\TO^{G_{t_2}}\right)\oplus 
H^1\left(\frac{G_{t_2}}{{G_{t_1}}},\TO^{G_{t_1}}\right)^{\frac{G_{t_3}}{G_{t_2}}}\oplus
\left(H^1(G_{t_1},\TO)^{\frac{G_{t_2}}{G_{t_1}}} \right)^{\frac{G_{t_3}}{G_{t_2}}}. 
\]
Using induction based on  (\ref{ind-33}) we obtain:
\[
H^1(G_1,\TO) \subseteq \bigoplus_{i=1}^f 
H^1\left(
\frac{G_{t_i}}{G_{t_{i-1}}},\TO^{G_{t_{i-1}}}
\right)^{\frac{G_{t_f}}{G_{t_i}}}\subseteq
 \bigoplus_{i=1}^f
 H^1\left(
\frac{G_{t_i}}{G_{t_{i-1}}},\TO^{G_{t_{i-1}}}
\right).
\]
and the desired result follows.
\end{proof}
Notice that in the above  proposition 
$G_{t_{i-1}}$  appears in the ramification filtration of $G_0$ thus 
the corollary to proposition IV.1.3 in \cite{SeL} implies that the ramification 
filtration of $G_{t_i}/G_{t_{i-1}}$ is constant. Namely, 
if $Q=G_{t_i}/G_{t_{i-1}}$ the ramification filtration of 
$Q$ is given by:
\[
Q_0=Q_1=\cdots=Q_{t_i}>\{1\}.
\]  Therefore,
$\delta_{t_1}=\dim_k H^1( G_n,\TO)$, and 
$\dim_k H^1(G_{t_i}/G_{t_{i-1}},\TO^{G_{t_{i-1}}})$  can be computed explicitly
by proposition \ref{comp-el-ab} since $G_n=G_{t_1}$, $G_{t_i}/G_{t_{i-1}}$, are
elementary abelian groups. Namely we will prove the following 
\begin{pro} \label{comp-quot-coho}
Let $\log_p(\cdot)$ denote the logarithmic function with base $p$.
Let $s(\lambda)=\log_p |G_{t_\lambda}|/|G_{t_{\lambda-1}}|$ and let $\mu_i$ be as in proposition 
\ref{comp-in-tan}. Then 
\[
\dim_k H^1
\left(
\frac{G_{t_\lambda}}{G_{t_{\lambda-1}}},
\TO^{G_{t_{\lambda-1}}}
\right)
=\sum_{i=1}^{s(\lambda)}
\left(
\lf
\frac{
(t_\lambda+1)(p-1)+a_i}{p}
\rf
-\lc \frac{a_i}{p} \rc
\right),
\]
 where $a_1=-t_\lambda-1+\mu_{\lambda-1}$, and $a_i= \lc \frac{a_{i-1}}{p} \rc$.
\end{pro}
\begin{proof}
The module $\TO^{G_{t_{\lambda-1}}}$ is computed in proposition \ref{comp-in-tan} to be   
isomorphic to $\pi_{\lambda-1}^{\mu_{\lambda-1}} \frac{d}{d\pi_{\lambda-1}}$,
 which in turn is $\frac{G_{t_\lambda}}{G_{t_{\lambda-1}}}$-equivariantly
isomorphic to $\pi_{\lambda-1}^{-t_\lambda-1+\mu_{\lambda-1}} k[[\pi_{\lambda-1}]]$. The desired result follows by using  proposition \ref{comp-el-ab}.
\end{proof}
\begin{rem}
If $n=1$, {\em i.e.} $G_2=\{1\}$ then the left hand side 
and  the right hand of (\ref{loc-def-dim}) are equal and 
the bound becomes the formula  in \cite{CK}. 
\end{rem}
\begin{pro}
We will follow the notation of \ref{comp-quot-coho}.
Suppose that for every $i$, $\frac{G_{t_i}}{G_{t_{i-1}}}$ is a cyclic $p$-group. Then the following equality 
holds:
\begin{eqnarray*}
\dim_k H^1(G_1,\TO) = & \sum_{i=1}^f \dim_k H^1(\frac{G_{t_i}}{G_{t_{i-1}}},\TO^{G_{t_{i-1}}})^{G_{t_f}/G_{t_i}} \leq \\
                    = &
\sum_{i=1}^f 
\left(
\lf
\frac{
(t_i+1)(p-1)-t_i-1+\mu_{i-1}}{p}
\rf
-\lc \frac{-t_i-1+\mu_{i-1}}{p} \rc
\right).
\end{eqnarray*}
\end{pro}
\begin{proof}
The kernel of the transgression at each step  is by lemma \ref{zero-trans} the whole 
$H^1(\frac{G_{t_i}}{G_{t_{i-1}}},\TO^{G_{t_{i-1}}})^{G_{t_f}/G_{t_i}}$. Therefore the right inner inequality in 
equation (\ref{loc-def-dim}) is achieved. The other inequality is trivial  by the computation done in 
proposition \ref{comp-quot-coho} but it is far from being best possible.
\end{proof}
%\begin{rem}
%All dimensions of cohomology groups that appear in the bound of 
%proposition \ref{loc-def-pro} can be made explicit with the aim 
%of proposition \ref{comp-in-tan} and by (\ref{eq-el-ab}). 
%\end{rem}
%
%
%
%
\section{Global Computations}\label{global-comp}
We consider the Galois cover of curves
$\pi: X \rightarrow Y=X/G$, and let $b_1,...,b_r$ be the 
ramification points of the cover.   
We will denote by
\[
e_0^{(\mu)} \geq e_1^{(\mu)} \geq e_2^{(\mu)}  \geq e_{n_\mu}^{(\mu)} > 1
\]
the orders of the higher ramification groups at the point $b_\mu$.  
The ramification divisor $D$ of the above cover  is a divisor 
supported at the ramification points $b_1,...,b_r$ and equals to 
\[
D=\sum_{\mu=1}^r \sum_{i=0}^{n_\mu} (e_i^{(\mu)}-1) b_\mu. 
\]
Let $\Omega_X^1$, $\Omega_Y^1$ be the sheaves of holomorphic differentials 
at $X$ and $Y$ respectively. The following formula holds 
\cite[IV. 2.3]{Hartshorne:77}:
\[
\Omega_X^1 \cong \O_X(D) \otimes \pi^* (\Omega_Y^1)
\]
and by taking duals
\[
\T_X\cong \O_X(-D)\otimes \pi^*(\T_Y)
\]
Thus $\pi_*(\T_X)\cong \T_Y \otimes \pi_*(\O_X(-D))$ and 
$\pi_*^G(\T_X)\cong \T_Y \otimes (\O_Y \cap \pi_*(\O_X(-D)))$. 
We compute (similarly with \cite[prop. 1.6]{CK}), 
\[
\pi_*^G(\T_X) =\T_Y \otimes \O_Y \left( -\sum_{\mu=1}^r 
\lc \sum_{i=0}^{n_\mu} \frac{(e_i^{(\mu)}-1)}{e_0^{(\mu)}} \rc b_i 
\right)
\]
Therefore, the global contribution to $H^1(G,\T_X)$ is given by 
\begin{align*}
H^1(Y,\pi_*^G(\T_X)) & \cong  H^1(Y, \T_Y \otimes \O_Y 
 \left( -\sum_{\mu=1}^r 
\lc \sum_{i=0}^{n_\mu} \frac{(e_i^{(\mu)}-1)}{e_0^{(\mu)}} \rc b_i
\right)
\\
&\cong
H^0(Y,\Omega_Y^{\otimes 2}  \left( -\sum_{\mu=1}^r 
\lc \sum_{i=0}^{n_\mu} \frac{(e_i^{(\mu)}-1)}{e_0^{(\mu)}} \rc b_i
\right)
\end{align*}
and by Riemann-Roch formula 
\begin{equation}\label{glob-contr}
\dim_k H^1(Y,\pi_*^G (\T_X)) = 3 g_Y -3 + \sum_{\mu=1}^r 
\lc \sum_{i=0}^{n_\mu} \frac{(e_i^{(\mu)}-1)}{e_0^{(\mu)}} \rc.
\end{equation}
On the other hand the local contribution can be  bounded by 
proposition \ref{loc-def-pro} and by combining the local and 
global contributions, we arrive at the desired  bound for the dimension.

\subsection{Examples}

Let $V=\Z/p\Z\oplus \cdots \oplus \Z/p\Z$ be an elementary abelian group acted on by the group $\Z/n\Z$. 
Assume that $G:=V \rtimes \Z/n\Z$ acts on the local field $k((t))$ and assume that the ramification filtration is 
given by 
\[
G_0>G_1=G_2=\cdots G_j > G_{j+1}=\{1\}.
\] 
Let $H:=\Z/p\Z$ be the first summand of $V$. Let $\sigma$ be a generator of the cyclic group $\Z/n\Z$ and assume that
$\sigma(t)=\zeta t$, where $\zeta$ is a primitive $n$-th root of one. 
The inflation-restriction sequence implies the short exact sequence:
\[
0 \rightarrow H^1(G/H ,x^a k[[x]]\frac{d}{dx}) \rightarrow H^1(G, t^{a'} k[[t]]\frac{d}{dt} ) 
\rightarrow H^1(H, t^{a'} k[[t]] \frac{d}{dt}) \rightarrow 0
\]
where $x^a k[[x]] \frac{d}{dx}=\left( t^{a'} k[[t]]\frac{d}{dt}\right)^H$.
The group $\Z/n\Z$ acts on $t^{a'}k[[t]]$ but there is no apriori well defined action of $\Z/n\Z$ on $x^a k[[x]]\frac{d}{dx}
=\left( t^{a'} k[[t]] \frac{d}{dt} \right)^H$, 
since the group $H$ might not be normal in $G$.
An element $d\in H^1(G/H, x^a k[[x]] \frac{d}{dx})$ is send by the inflation map on the $1$-cocycle 
$\mathrm{inf}(d)$ that is a map
\[\mathrm{inf}(d): G \rightarrow  t^{a'} g(t) \frac{d}{dt} \in t^{a'} k[[t]] \frac{d}{dt}, \] 
and  the action of $\sigma$  can be considered on the image of the inflation map, 
sending $\mathrm{inf}(d)(g) \mapsto \sigma(\mathrm{inf}(d(g))$.
We observe that $\sigma(\mathrm{inf})(d(g)))$ is zero for any $g\in H$, by the definition of the inflation map, therefore there 
is an element $a\in t^{a'} k[[t]] \frac{d}{dt}$ such that 
\[
\sigma(\mathrm{inf}{d}(g))+ a^g -a \in x^a k[[x]] \frac{d}{dx},
\]
therefore  we can consider the element 
\[
\sigma(\mathrm{inf}(d))+ a^g -a= \mathrm{inf}(d').
\]
This means that although  there is no well defined action of $\Z/n\Z$ on $k[[x]]$ we can define 
 $\sigma(d)=d' $ modulo cocycles. 
In what follows we will try to compute the element $d' \in H^1(G/H, x^a k[[x]] \frac{d}{dx})$. 

Assume that the Artin-Schreier extension $k((t))/k((x))$ is given by the equation $1/y^p -1/y =1/x^j$. Then, 
we have computed that if $g$ is a generator if $H$ then  
\[
g(t)= \frac{t}{(1+t^{j})^{1/j}}, \mbox{ and } x=\frac{t^p}{(1-t^{j(p-1)})^{1/j}}.
\] 
The action of $\sigma$ on $x$, where $x$ is seen as an element in $k[[t]]$  is given by 
\[
\sigma(x)=\sigma \frac{t^p}{(1-t^{j(p-1)})^{1/j}}= \zeta^p x \frac{ (1-t^{j(p-1)})^{1/j}}{(1-\zeta^{j(p-1)}t^{j(p-1)})^{1/j}}=\zeta^p x u,
\]
where $u=\frac{ (1-t^{j(p-1)})^{1/j}}{(1-\zeta^{j(p-1)}t^{j(p-1)})^{1/j}}$ is a unit of the form $1+y$, where $y\in t^{2j(p-1)}k[[t]]$.

The cohomology group $H^1(V/H, x^a k[[x]])$ is generated by the elements $\{1/x^\mu, b\leq \mu \leq a, \binom{ \mu/n}{p-1}=0\}$. Each element
$1/x^\mu$ is written as $1/x^{\mu} x^{j+1} \frac{d}{dx}$ and it is lifted to
\[
1/x^{\mu} t^{j+1} \frac{d}{dt}\stackrel{\sigma}{\longrightarrow} \zeta^{-p\mu +j} 1/x^{\mu} u^{-\mu} t^{j+1} \frac{d}{dt}.
\]
In the above formula we have used the fact that the adjoint action of 
$\sigma$ on 
$t^r\frac{d}{dt}$ is given by 
$t^r\frac{d}{dt}\stackrel{\sigma}{\rightarrow} \zeta^{(r-1)}t^r\frac{d}{dt}$
\cite[3.7]{CK}.

Obviously the unit $u$ is not $H$-invariant but we can add to u a $1$-coboundary so that it becomes the $H$-invariant element $\mathrm{inf}(d')$.
We observe that this coboundary is of the form $a^g-a$, and obviously $a^g-a$ has to be in $t^{2 {p-1} j}k[[t]]$.
This gives us that  
\[
\left(1/x^\mu\right)'  = \zeta ^{\mu} 1/x^\mu + o,
\]
where $o$ is a sum of terms $1/x^\nu$ with $-a<\nu$  and therefore $o$ is cohomologous to zero.
Using induction one can prove the following 
\begin{lemma} \label{1-action-zeta}
Let $1/\pi_i^{i_\lambda}$, $\lambda=1,\ldots,s$, $b_i \leq i_\lambda \leq -a_i$ so that $\binom{i_\lambda/n}{p-1}=0$ and 
$b_i=-a_i-j$ if $p\mid a_i$, $b_i=-a_i-j+1$ if $p \nmid a_i$ be the basis elements of the cohomology group $H^1(V,\TO$. 
Then the action of the generator $\sigma \in \Z/n\Z$ on $\TO$ induces the following action on the basis elements:
\[
\sigma(\frac{1}{\pi_i^\mu})=\zeta^{-{p^i}\mu+j} \frac{1}{\pi^\mu}.
\]
\end{lemma}

$\;$\\
{\bf 1.}
The Fermat curve 
\[F:x_0^n +x_1^n +x_2^n=0\] defined over an 
algebraically closed field $k$ of characteristic $p$, 
such that $n-1=p^a$ is a power of the characteristic 
is a very special curve. 
Concerning its  automorphism group, the Fermat curve 
has maximal automorphism group with respect to the genus \cite{StiI}. 
Also it leads to Hermitian function fields, that are optimal 
with  respect to the number of $\F_{p^{2a}}$-rational points
and Weil's bound. 

It is known that the Fermat curve is totally supersingular, 
{\em i.e.}, the Jacobian variety $J(F)$ of $F$ has $p$-rank 
zero, so this curve cannot be studied by the tools of \cite{CK}.
The group of automorphism of $F$ was computed in 
\cite{Leopoldt:96} to be the projective unitary group 
$G=PGU(3,q^2)$, where $q=p^a=n-1$. H.Stichtenoth \cite[p. 535]{StiI}
proved that in the extension $F/F^G$ there are two ramified points 
$P$,$Q$ and one is wildly ramified and the other is tamely ramified.  
For the ramification group $G(P)$ of the wild ramified point $P$
we have that $G(P)$ consists of the $3\times 3 $ matrices of the form
\begin{equation}
\label{paparia}
\begin{pmatrix}
1 & 0 & 0 \\
\alpha & \chi & 0 \\
\gamma & -\chi \alpha^q & \chi^{1+q}
\end{pmatrix},
\end{equation}
where $\chi,\alpha,\gamma\in \F_{q^2}$ and 
$\gamma+\gamma^q=\chi^{1+q}-1-\alpha^{1+q}$. Moreover Leopoldt proves that 
the order of $G(P)$ is $q^3(q^2-1)$ and 
the ramification filtration is given by 
\[
G_0(P) > G_1(P) >G_2(P) =\cdots = G_{1+q}(P) > \{1\},
\]
where
\[G_1(P)=\ker(\chi:G_0(P) \rightarrow \F_{q^2}^*)\]
and
\[
G_2(P)=\ker(\alpha:G_1(P) \rightarrow \F_{q^2}).
\]
In this section we will compute the dimension  of tangent space of the 
global         deformation functor. 
Namely, we will prove:
\begin{pro}
Let $p$       be a prime number, $p>3$ let $X$ be  the Fermat curve 
\[
x_0^{1+p} + x_1^{1+p} +  x_2^{1+p}=0.
\]  Then $\dim_k H^1(X,G,\T_X)=0$. 
\end{pro}
\begin{proof}
By the assumption $q=p$ and by the computations of Leopoldt 
mentioned above we have $G_2=\cdots =G_{p+1}=\Z/p\Z$. The different of 
$G_{p+1}$ is computed $(p+2)(p-1)$. Hence, 
according to  (\ref{cyc-coh}) 
\[
\dim_k H^1(G_{p+1},\TO)= \lf \frac{(p+2)(p-1)-(p+2)}{p} \rf
-\lc  \frac{-p-2}{p} \rc= p
\]
 Proposition \ref{final-cyc-comp} implies that the set
\[
\left\{
\frac{1}{t^i},   2 \leq i \leq p+2 \mbox{ where }  \binom{ \frac{i}{p+1}}{p-1}=0
\right\}
\]
is a $k$-basis of $H^1(G_{p+1},\TO)$. Indeed, the group $G_{1+p}$ has conductor $1+p$ and 
 $\TO$ is $G_{1+p}$-equivariantly isomorphic to $t^{-p-2}k[[t]]$. Thus following the 
 notation of proposition \ref{final-cyc-comp} $-a=p+2$ and $b=2$.
 The rational number $(1+p)^{-1}$ has the following $p$-adic expansion:
 \[
 \frac{1}{1+p} =1+ (p-1) p + (p-1) p^3 + (p-1) p^5 + \ldots
 \]
 and using lemma \ref{binom-zero} we obtain that for $2\leq i \leq p+2$ the only integer $i$  such that 
   $\binom{ \frac{i}{p+1}}{p-1}\neq 0$ is $i=p-1$. Thus, the elements
  \[
\left\{
\frac{1}{t^i},   2 \leq i \leq p+2, i\neq p-1
\right\}
\]
 form a $k$-basis of $H^1(G_{p+1},\TO)$.
 
   Leopoldt in \cite{Leopoldt:96} proves that the $G_0(P)$ acts faithfully on the $k$-vector space
   $L((p+1)P)$ that is of dimension $3$ with basis functions $1,v,w$ and the representation 
   matrix is given by (\ref{paparia}). Moreover, 
   the above functions have $t$-expansions of the following form 
   $v=\frac{1}{t^p} u$, where $u$ is a unit in $k[[t]]$ and $w=\frac{1}{t^{p+1}}$, 
   for a suitable choice of the local uniformizer $t$ at the point $P$. 
The functions $v,w$ generate the function field coresponding to the Fermat curve and they satisfy
the relation $u^n=w^n-(w+1)^n$, therefore one can compute that the unit $u$ can be written as 
\[
u=1+ t^{p+1}g, \;\; g\in k[[t]].
\]
Let $\sigma$ be an element given by a matrix as in equation (\ref{paparia}).
   The action of $\sigma\in G_1=G_1(P)$ on powers of $\frac{1}{t}$ is given by 
   \begin{equation} \label{G1action}
   \frac{1}{t^i} =\frac{(1+ \gamma t^{p+1} - a ^q u t)^\frac{i}{p+1}}{t^i}. 
   \end{equation}
and the action on the basis elements $\{1/t^i,\; 2\leq i \leq p+2, i\neq p-1\}$ is given by
\[
\frac{1}{t^i} \mapsto \frac{1}{t^i}+ \sum_{\nu=1}^{i-2} a^{q\nu} \binom{\frac{i}{p+1}}{\nu} \frac{1}{t^{i-\nu}}.
\]
We observe that the matrix of this action  is given by 
\[
A_\sigma=
% use packages: array
\left(
\begin{array}{lllll}
1 & 0 & 0 & 0 & 0 \\ 
\frac{2}{p+1} & 1 & 0 & 0 & 0 \\ 
 & \frac{3}{p+1} & 1 & 0 & 0 \\ 
* & * & \ddots & 1 & 0 \\ 
* & * & * & \frac{p+2}{p+1} & 1
\end{array}
\right)
\]
We observe that $\sigma(1/t^2)=1/t^2$ and $\sigma(1/t^p)=1/t^p$, and moreover that all elements
below the diagonal of the matrix $A_\sigma$ are $\frac{i}{p+1}$ and are non-zero unless $i=p$.
Therefore the eigenspace of the eigenvalue $1$ is 2-dimensional, and we can give a basis:
\[
H^1(G_{1+p},\TO)^{G_1/G_{1+p}}=_k\left\langle  \frac{1}{t^2},\frac{1}{t^p}\right\rangle
\]
In order to compute $H^1(G_1(P),\TO)$ we consider the exact sequence
\[
1 \rightarrow G_2 \rightarrow G_1 \stackrel{\alpha}{\longrightarrow}
G_1/G_2\cong \Z/p\Z \times \Z/p\Z \rightarrow 1
\]
and the corresponding low-degree-term Lyndon-Hochschield-Serre sequence. 
The  group $G_2$ is of conductor $p+1$ thus $\TO^{G_2}=\TO^{\Z/p\Z}$ is given by proposition \ref{comp-in-tan} ($p>2$):
\[
\TO^{G_2}=x^{p+2-\lf \frac{p+2}{p} \rf} k[[x]]\frac{d}{dx}=x^{p+1}k[[x]]\frac{d}{dx},
\] where $x$ is a 
local uniformizer for $\O^{G_2}$.
By \cite[Cor. p.64]{SeL} the 
ramification filtration for  $G_2/G_1$ is 
\[
G_0/G_2 > G_1/G_2  > \{1\},
\]
hence the different for the subgroup $\Z/p\Z$ of $G_2/G_1$ is $2(p-1)$, and the conductor 
equals $1$.
Lemma \ref{TO-invariants2} implies
 $x^{p+1}k[[x]] \frac{d}{dx}$ is $G_1/G_2$-equivariantly isomorphic to $x^{p+1-2}k[[[x]]$.
Therefore,
\[H^1(G_1/G_2,\TO^{G_2})=H^1(\Z/p\Z , x^{p-1} k[[x]]) \oplus 
                         H^1(\Z/p\Z, (x^{p-1} k[[x]])^{\Z/p\Z} )
			 \]

We compute
\[
\dim_k H^1(\Z/p\Z, x^{p-1} k[[x]])=\lf \frac{2(p-1)+ p-1}{p}\rf -\lc 
\frac{p-1}{p}\rc=1. 
\]
On the other hand, if $\pi$ is a local uniformizer for $k((x))^{\Z/p\Z}$
then 
\[
(x^{p-1} k[[x]])^{\Z/p\Z}=\pi^{\lc  \frac{p-1}{p} \rc} k[[\pi]] =\pi
k[[\pi]] 
\]
and the dimension of the cohomology group is computed:
\begin{align*}
\dim_k H^1(\Z/p\Z, (x^{p-1} k[[x]])^{\Z/p\Z} ) 
&=\dim_k H^1(\Z/p\Z,\pi k[[\pi]])=
\\
&=
\lf \frac{2(p-1)+1}{p} \rf -\lc \frac{1}{p}\rc = 0.
\end{align*}
Using the bound for the kernel of the transgression we see that 
\begin{eqnarray} \nonumber
1=\dim_k H^1(\frac{G_1}{G_2},\TO^{G_2})
\leq \dim_k H^1(G_1 ,\TO) \leq  & 
\dim_k H^1(\frac{G_1}{G_2},\TO^{G_2})+ \\ +         & \dim_k H^1(G_2,\TO)^{\frac{G_1}{G_2}}
=3. \label{lastbound}
\end{eqnarray}

In order to compute the action of $G_0$ on $G_1/G_2$ we observe that 
\begin{equation} \label{conj-fin12}
\begin{pmatrix}
1 & 0 & 0 \\
a & \chi  & 0 \\
* & -\chi a^p \\
\end{pmatrix}
\begin{pmatrix}
1 & 0 & 0 \\
b & 1 & 0 \\
* & -b^p & 1 \\
\end{pmatrix}
\begin{pmatrix}
1 & 0 & 0 \\
a & \chi  & 0 \\
* & -\chi a^p \\
\end{pmatrix}=
\begin{pmatrix}
1 & 0 & 0 \\
\chi b & 1 & 0 \\
* & - \chi^{q+1} b^p & 1 \\
\end{pmatrix}
\end{equation}
If $\begin{pmatrix} 1 & 0 & 0 \\
b & 1 & 0 \\
* & -b^p & 1 \\
\end{pmatrix}
$ is an element of $\Z/p\Z\cong \F_p \subset \F_{p^2}$ then $b^p=b$. By looking at the computation of (\ref{conj-fin12}) we 
see that  the conjugation action of $G_0/G_1$ to $\F_p$ is given 
by multiplication 
$b\mapsto \chi^{1+q} b$. Observe that $(\chi^{1+p})^{p-1}=\chi^{p^2-1}=1$, thus $\chi^{1+p}\in \F_p$. 
The action on the cocycles is given by sending the cocycle 
$d(\tau)$
to $d(\sigma \tau \sigma^{-1})^{\sigma^{-1}}$ therefore the basis cocycle $\frac{1}{x^{1-p}}$ of the one dimensional cohomology group 
$H^1(\Z/p\Z, x^{p-1} k[[x]])$   
$d$ goes to $\chi^{p(p-1)+1+1+p} d=\chi^{-p^2} d$  under the conjugation action, as 
one sees by applying lemma \ref{1-action-zeta}.  Lemma \ref{tame-computing} implies that $H^1(G_1/G_2,\TO^{G_2})^{G_0/G_1}=0$.
Similarly the conjugation action of $G_0/G_1$ on an element of $G_2$ can be computed to be multiplication 
of $\tau$ by $\chi^{1+p}\in \mathbb{F}_p$, and the same argument shows that ${H^1(G_{1+p},\TO)^{G_1/G_{1+p}}}^{G_0/G_1}=0$.

Finally the global contribution is computed  by formula (\ref{glob-contr})
\begin{align*}
\dim_k H^1 (F^G,\pi_*(\T_F)) &=-3+ 
\lc \sum_{i=0}^{p+2} \frac{|G(P)_i|-1}{|G(P)|} \rc  +
\lc 1-\frac{1}{|G(Q)|}  \rc=
\\
&=-3 + 2 +1=0
\end{align*}
The fact that 
the tangent space of the deformation functor is zero dimensional 
is compatible with  the fact that
there is only one isomorphism class of curves  $C$ such that 
$|Aut(C)| \geq 16 g_C^4$  \cite{StiI}. 
\end{proof}
{\bf 2. $p$-covers of $\mathbb{P}^1(k)$}
We consider curves $C_f$ of the form  
\[
C_f:w^p-w= f(x),
\]
where $f(x)$ is a polynomial of degree $m$. We will say that such a curve is in reduced 
form if the polynomial $f(x)$ is of the form 
\[
f(x)=\sum_{i=1,(i,p)=1}^{m-1} a_i x^i+x^m.
\]
Two such curves $C_f$, $C_g$ in reduced form are isomorphic if and only if 
$f=g$. 
 The group $G:=Gal\big(C_f/\mathbb{P}^1(k)\big)\cong \Z/p\Z$ acts on $C_f$. 
 We observe that the number of independent monomials $\neq x^m$ in the above sums is given by:
 \begin{equation} \label{pcovers1}
 m - \lf \frac{m}{p} \rf-1, 
 \end{equation}
 since $\#\{ 1\leq i \leq m, p\mid i \}= \lf \frac{m}{p} \rf $.
 
 We will compute the tangent space of the  deformation functor of   the curve 
 $C_f$ together with the group $C_f$. 
Let $P$ be the point above $\infty \in \mathbb{P}^1(k)$. This is the only point 
that ramifies in the cover $C_f \rightarrow \mathbb{P}^1(k)$, and the group $G$
admits the following ramification filtration:
\[
G_0=G_1=G_2=\cdots=G_m > G_{m+1}=\{1\}.
\]
The different is computed $(p-1)(m+1)$ and $\TO\cong t^{-m-1}k[[t]]$. 
Thus the space $H^1(G,\TO)$ has dimension $d$ 
\[
d=\lf 
 \frac{(p-1)(m+1)-(m+1)}{p}
 \rf -
\lc \frac{-(m+1)}{p} \rc =m+1- \lc \frac{2m+2}{p} \rc + \lf \frac{m+1}{p} \rf.
\] 
Let $a_0 + a_1 p + a_2 p^2 +\cdots$
be the $p$-adic expansion of $m+1$. We observe that 
\[
\lc \frac{2m+2}{p} \rc - \lf \frac{m+1}{p}\rf= \lc \frac{2a _0}{p} + \sum_{i\geq 1} 2 a_i p^{i-1} \rc - \sum_{i\geq 1} a_i p^{i-1},  
\] 
therefore, if $p\nmid m+1$
\[
\lc \frac{2m+2}{p} \rc - \lf \frac{m+1}{p}\rf= \lf \frac{m+1}{p} \rf + \delta, 
\]
where 
\[
\delta=
\left\{
\begin{array}{ll}
2 & \mbox{ if } 2a_0> p\\
1 & \mbox { if } 2a_0 <p.
\end{array}
\right.
\]
Thus, 
\[
d=\left\{
\begin{array}{ll}
m+1 - \lf \frac{m+1}{p}\rf & \mbox{ if } p\mid m+1 \\
m-  \lf \frac{m+1}{p}\rf -\delta &\mbox{ otherwise}
\end{array}
\right.
\]
Finally, we compute that 
\[
dim_k H^1(Y,\pi_*^G (\T_X)) = -3 +  
\lc \frac{(m+1)(p-1)}{p} \rc=m-2 - \lf \frac{m+1}{p} \rf.
\]

{\bf 3. Lehr-Matignon Curves.}
Let us consider the curve
\[
C:y^p -y = \sum_{i=0}^{m-1} t_i x^{1+p^i} + x^{1+p^m} ,
\]
defined over the algebraically closed field $k$ of characteristic $p >2$. Let $n=1+p^m$ denote the 
degree of right hand side of the above equation.
The automorphism group of these curves were studied by Matignon-Lehr in \cite{CL-MM} and these curves  were 
considered also by G. van der Geer, van der Vlught in \cite{Geer-Vlugt92} in connection with coding theory.
Notice that the extreme Fermat curves studied in example 1 can be written in this form by a suitable transformation \cite[VI.4.3 p.203]{StiBo}.
Let $H=Gal(C/\mathbb{P}^1(k))$. The automorphism group $G$  of $C$ can be expressed in the form 
\[
1 \rightarrow H \rightarrow G \rightarrow V \rightarrow 1,
\]
where $V$ is the vector space of roots of the additive polynomial 
\begin{equation} \label{adjMat}
\sum_{0\leq i \leq m} \big(t_i^{p^{m-i}}Y^{p^{m-i}}+ t_i^{p^m}Y^{p^{m+i}}\big)
\end{equation} \cite[prop. 4.15]{CL-MM}.
Moreover there is only one point $P\in C$  that ramifies in the cover  $C\rightarrow C^G$, namely the point above $\infty \in \mathbb{P}^1(k)$.

In order to simplify the calculations 
we assume that $t_0=\cdots=t_{m-1}=0$ so the curve is given by
\begin{equation} \label{defsimple}
y^p-y= x^{p^m+1}.
\end{equation}
The polynomial in (\ref{adjMat})  is given by $Y^{p^{2m}}+Y$ and the vector space $V$ of the 
roots is $2m$-dimensional.
Moreover, according to \cite{CL-MM} any automorphism $\sigma_v$ coresponding to  $v\in V$ is given by
\[
\sigma_v(x)=x+v, \;\;\;\sigma(y)=y + \sum_{\kappa=0}^{m-1} v^{p^{m+k}}x^{p^k}. 
\] 
Observe that $w$ (resp. $x$)  has a unique pole of order $p^m+1$ (resp. $p$) at the point above $\infty$, 
 so we can select the local 
uniformizer $\pi$ so that 
\[
y=\frac{1}{\pi^{p^m+1}}, \;\;\; x=\frac{1}{\pi} u,
\]
where $u$ is a unit in $k[[\pi]]$. By replacing $x,y$ in (\ref{defsimple}) 
we observe that the unit $u$ is of the form $u=1+ \pi^{p^m}$.

A simple computation  based on the basis $\{1,x,\ldots, x^{p^{m-1}},y\}$, of the vector space 
 $L\big( (1+p^m)P\big)$ given in \cite[prop. 3.3]{CL-MM}.
 shows that the ramification filtration of $G$ is 
\[
G=G_0=G_1 > G_2=\ldots = G_{p^m+1}  > \{1\},
\]
where $G_2=H$ and $G_1/G_2=V$. 
Using proposition \ref{final-cyc-comp} we obtain the following basis for $H^1(G_2,\TO)$:
\[
\left\{
\frac{1}{\pi_1^i},\; 2 \leq i \leq p^m+2, \mbox{such that} \binom{\frac{i}{p^m+1}}{p-1}=0.
\right\}.
\]
We have to  study the action of $G_1/G_2$ on $H^1(G_2,\TO)$. 
From the action of $\sigma_v$ on $y$ we obtain that that the action on the basis elements of 
$H^1(G_{p^m+1},\TO)$ is given by 
\[
\sigma_v(\frac{1}{\pi^i})=\frac{1}{\pi^i} + 
\frac{
\left(
\sum_{\kappa=0}^{m-1} u^{p^k} v^{p^{m+\kappa}} \pi^{p^m+1-p^{k+1}}
\right)^{\frac{i}{p^m+1}}
}{\pi^i}=
\]
\[
\frac{1}{\pi^i} + \frac{i}{p^m+1} v^{p^{2m-1}} \frac{1}{\pi^{i-1}} + \ldots.
\]
We observe that if $p\mid i$ then all binomial coefficients $\binom{\frac{i}{p^m+1}}{\kappa}$ that 
contribute  a coefficient $1/\pi^\kappa$, $2\leq \kappa \leq p^m+2$ are zero.
Therefore, the elements $1/\pi^2,1/\pi^{\nu p}$ are invariant.
Moreover, by writing down the action of $\sigma_v$ as a matrix we see that there are 
no other invariant elements, so the dimension is computed $(p>2)$:
\[
\dim_k H^1(G_{p^m+1},\TO)^{G_1/G_{p^m+1}}=1 + \lf \frac{ p^m+2}{p}\rf=1 + p^{m-1}.
\]
We observe that this dimension coincides with the computation done on the Fermat curves $m=1$.

We proceed by computing $H^1(V,\TO^H)$. The space $\TO^H$ is computed by proposition  \ref{comp-el-ab}
\[x^{p^m+2 -\lf \frac{p^m+2}{p}\rf}k[[x]]\frac{d}{dx}=x^{p^m+2-p^{m-1}}k[[x]]\frac{d}{dx}.\]
Thus,
\[
\dim_k H^1(V,\TO^H)=\sum_{\nu=1}^{2m} \left(
\lf 
\frac{2(p-1)+a_i}{p}
\rf
-
\lc
\frac{a_i}{p}
\rc
\right),
\]
where $a_1=p^m-p^{m-1}$, and $a_i=\lc \frac{a_{i-1}}{p} \rc$.
By computation $a_\nu=p^{m-\nu+1}-p^{m-\nu}$ for  $1\leq \nu \leq m$, and $a_\nu=1$ for $\nu>m$. 
On the other hand, an easy computation shows that 
\[
\lf 
\frac{2(p-1)+a_i}{p}
\rf
-
\lc
\frac{a_i}{p}
\rc
= 
\left\{
\begin{array}{ll}
1 & \mbox{ if } 1\leq \nu < m \\
2 & \mbox{ if } \nu=m\\
0 & \mbox{ if } m < \nu
\end{array}
\right.,
\]
thus the dimension of the tangent space is $m+1$. 

We have proved that the dimension of $H^1(G_1,\TO)$ is bounded by 
\[
m+1=\dim_k H^1(G_1/G_2,\TO^{G_2}) \leq H^1(G_1,\TO) \leq 
\] 
\[
\dim_k H^1(G_1/G_2,\TO^{G_2})+  H^1(G_2,\TO)^{G_1/G_2}=2+m+p^{m-1}.
\]
Unfortunately we can not be more precise here: an exact computation involves the computation of the kernel of the transgression and 
such a computation requires new ideas and tools.

To this dimension we must add the contribution of  
\[
\dim_k H^1(Y,\pi_*^G (\T_X)) = 3 g_Y -3 + \sum_{\kappa=1}^r 
\lc \sum_{i=0}^{n_\kappa} \frac{(e_i^{(\kappa)}-1)}{e_0^{(\kappa)}} \rc=
\]
\[
-3+ \lc 2\frac{p^{m+1}-1}{p^{m+1}} + m\frac{p^m-1}{p^{m+1}} \rc= 
-1 + \lc \frac{m}{p} -\frac{2+m}{p^{m+1}} \rc.
\]
The later contribution is $>0$ if $m\gg p$.

{\bf 4. Elementary Abelian extensions of $\mathbb{P}^1(k)$.}\label{examp4} 
Consider the curve $C$ so that $G_0=(\Z/p\Z)^s \rtimes \Z/n\Z$ is the ramification group of wild ramified point,
and moreover the ramification filtration is given by 
\begin{equation} \label{113456}
G_0 > G_1=\ldots= G_j > G_{j+1}=\{1\}.
\end{equation}
An example of such a curve is provided by the curve defined by
\begin{equation} \label{psn}
C: \sum_{i=0}^s a_i y^{p^i}=f(x), 
\end{equation}
where $f$ is a polynomial of degree $j$ and all monomial summands  $a_k x^k $ of $f$ have exponent 
congruent to $j$ modulo $n$. Let $V$ be the $\mathbb{F}_p$-vector space of the
roots of the additive polynomial $\sum_{i=0}^ s a_i y^{p^i}$.  Assume that the automorphism group 
of the curve defined by (\ref{psn}) is $G:=V \rtimes \mathbb{Z}/ n\Z$. Thus $C\rightarrow \mathbb{P}^1(k)$ is 
Galois cover ramified only above $\infty$, with ramification group $G$ and ramification filtration is computed to be as 
in equation (\ref{113456}).

Let us now return to the general case.
Let us denote by $V$ the group $(\Z/p\Z)^s$. The group $V$ admits the structure of a $\mathbb{F}_p$ vector space, where $\mathbb{F}_p$ is 
the finite field with $p$-elements.
The conjugation action of $\Z/n\Z$ on $V$ implies a representation 
\[
\rho: \Z/n\Z \rightarrow GL(V).
\]
Since $(n,p)=1$, Mascke's Theorem gives that $V$ is the direct sum of simple $\Z/n\Z$-modules, 
{\em i.e.}, $V=\oplus_{i=1}^r V_i$. 
On the other hand, lemma \ref{action-on-cocycles} implies that the conjugation action is given 
by multiplication by $\zeta^j$, where $\zeta$ is an apropriate primitime $n$-th  root of one and 
$j$ is the conductor of the extension. 
If $\zeta^j \in \mathbb{F}_p$ then all the $V_i$ are one dimensional. In the more general case one has 
to consider representations 
\[
\rho_i : \Z/n\Z \rightarrow GL(V_i),
\]
where $\dim_{\mathbb{F}_p}V_i=d$. Notice that the dimension $d$ is the degree of the extension 
$\mathbb{F}_q/\mathbb{F}_p$, where $\mathbb{F}_q$ is the smallest field containing  $\zeta^j$.
 Let $e_1^{(i)},\ldots,e_d^{(i)}$ be an $\mathbb{F}_p$-basis of $V_i$,
and let us denote by $(a_{\mu\nu}^{(i)})$ the entries of the matrix  corresponding to $\rho_i(\sigma)$, 
where $\sigma$ is a generator of $\Z/n\Z$. The conjugation action on the arbitrary 
\begin{equation}\label{per-not1}
v=\sum_{i=1}^r\sum_{\mu=1}^d \lambda^{(i)}_\mu e_{\mu}^{(i)} \in V
\end{equation}
 is given by:
\begin{equation} \label{conj-acti-fin11}
\sigma v \sigma^{-1}=  
\sum_{i=1}^r \sum_{\mu=1}^d  \lambda_\mu^{(i)}  \sigma e_\mu^{(i)} \sigma^{-1}= 
\sum_{i=1}^r  \sum_{\mu=1}^d \lambda_\mu^{(i)} \sum_{\nu=1}^d a_{\mu\nu}^{(i)}e^{(i)}_\nu.
\end{equation}
In particular, 
\begin{equation} \label{conj-acti-fin12}
\sigma e_\mu^{(i)} \sigma^{-1} = \sum_{\nu=1}^d a_{\mu\nu}^{(i)} e_{\nu}^{(i)}.
\end{equation}
For the computation of $H^1(G,\TO)$, we  notice first that the  group  $H^1(V,\TO)$ can be computed using 
proposition \ref{triv11} and the isomorphism $\TO\cong t^{-j-1}k[[t]]$.

Next  we consider the conjugation action of $\Z/n\Z$ on $H^1(V,\TO)$, in order to compute 
$H^1(G,\TO)=H^1(V,\TO)^{\Z/n\Z}$.  By (\ref{gen-fo}) we have 
\begin{equation} \label{decom-cocycl-123}
H^1(V,\TO)=\oplus_{\lambda=1}^s H^1(\Z/p\Z,\TO^{\oplus_{\nu=1}^{\lambda-1} \Z/p\Z}),
\end{equation}
{\em i.e.}, the arbitrary cocycle $d$ representing  a cohomolgy class in $H^1(V,\TO)$ can be writen 
as a  sum of cocycles $d_\nu$  representing cohomology classes in $H^1(\Z/p\Z,\TO^{\oplus_{\nu=1}^{i-1} \Z/p\Z})$.
Let us follow a similar to (\ref{per-not1}) notation for the decomposition of $d$,
and  write $d=\sum_{i=1}^r \sum_{\nu=1}^d b^{(i)}_{\nu} d^{(i)}_{\nu}$, where
 $d_\nu^{(i)} (e^{(j)}_\mu)=0$ if $i\neq j$ or $\nu\neq \mu$.
Therefore,
\begin{equation} \label{ton-poulo1}
d(\sigma e_{\mu}^{(i)} \sigma^{-1})=
d
 \left( \sum_{\nu=1}^d  a_{\mu\nu}^{(i)} e_{\nu}^{(i)}\right)=
 \sum_{\nu=1}^d b_\nu^{(i)} a_{\mu\nu}^{(i)} d_{\nu}^{(i)} (e_{\nu}^{(i)}).
\end{equation}
We have now  to compute the $\Z/n\Z$-action on $d^{(i)}_k$. 
By lemma \ref{1-action-zeta} the element $\sigma$ acts on the  basis elements  
$\frac{1}{\pi_i^\mu}$ of $H^1(V,\TO)$ as follows
\begin{equation} \label{1234}
\sigma(\frac{1}{\pi_i^\mu})=\zeta^{-{p^i}\mu+j} \frac{1}{\pi^\mu}.
\end{equation}

By the above remarks we arrive at 
\begin{equation}\label{ton-poulo2}
\sigma(d)(e_{\mu}^{(i)}):=d(\sigma e_{\mu}^{(i)} \sigma^{-1})^{\sigma^{-1}}=
\sum_{\nu=1}^d b_\nu^{(i)} a_{\mu\nu}^{(i)} \zeta^{-c(\nu,i)} d_\nu^{(i)} (e_{\nu}^{(i)}),
\end{equation}
where $c(\nu,i)$ is the appropriate exponent, defined in (\ref{1234}). Let us denote by 
$A^{(i)}$ the $d\times d$ matrix $(a_{\nu\mu}^{(i)})$. 
By (\ref{ton-poulo2}) $\sigma(d)(e_\mu^{(i)})=d(e_\mu^{(i)})$ if and only if
$b:=(b_1^{(i)},\ldots,b_d^{(i)})$ is a solution of the linear system 
\[
(A^{(i)}\cdot \mathrm{diag}^{-1}(\zeta^{c(1,i)},\zeta^{c(2,i)},\ldots,\zeta^{c(d,i)}) -\mathbb{I}_d)b=0.
\]
This proves that the dimension of the solution space is equal to the dimension of the eigenspace of the 
eigenvalue $1$ of the matrix: $A^{(i)}\cdot \mathrm{diag}^{-1}\big(\zeta^{c(1,i)},\zeta^{c(2,i)},\ldots,\zeta^{c(d,i)}\big)$. 

Moreover using a basis of the form $1,\zeta,\zeta^2,\ldots,\zeta^{d-1}$ for the simple space $V(i)$,
we obtain that
\[
A^{(i)}=
\left(
\begin{array}{ccccc}
0 & 0 & \cdots & 0 & -a_0 \\ 
1 & 0 &  &  &           -a_1 \\ 
0 & 1 & \ddots &  &   \vdots \\ 
0 &  &  1& 0 & -a_{d-1} \\ 
0 &  &  & 1 &  - a_d
\end{array}
\right)
\]
It can be proved by induction that the characteristic polynomial of $A^{(i)}$ is $x^{d}+\sum_{\nu=0}^{d-1} a_\nu x^\nu$, and 
under an appropriate basis change
$A^{(i)}$ can be a written in the form $\mathrm{diag}(\zeta^{j},\zeta^{2j},\ldots,\zeta^{j(d-1)})$. 
Moreover, the characteristic polynomial of the matrix 
$A^{(i)}\cdot \mathrm{diag}^{-1}\big(\zeta^{c(1,i)},\zeta^{c(2,i)},\ldots,\zeta^{c(d,i)}\big)$
can be computed inductively to be 
\[
f_i(x):=
x^d+\zeta^{c(d,i)} a_{d-1} x^{d-1} +\zeta^{c(d,i)+c(d-1,i)}a_{d-1}+
 \cdots+ 
  \zeta^{\sum_{\nu=2}^d c(\nu,i)} x+ 
  \zeta^{\sum_{\nu=1}^d c(\nu,i)} a_0
\]
If, $f_i(1)\neq 0$, then we set $\delta(i)=0$. If $f_i(1)=0$  we set $\delta(i)$ to be the 
multiplicity of the root $1$. The total invariant space has dimension
\[
\dim_k H^1(G,\TO) = \sum_{i} \delta(i).
\]

{\bf Comparison with the work of Cornelissen-Kato}
We will apply the previous calculation to the case of ordinary curves $j=1$ and we will obtain the formulas 
in \cite{CK}. We will follow the notation of proposition \ref{triv11}. The number $a_1=-j-1=-2$. Thus, 
$a_2=\lc -2/p \rc=-\lf 2/p \rf=0$ (recall that we have assumed that $p\geq 5$. Furthermore $a_i=0$ for $i\geq 2$.
For the numbers $b_i$ we have $b_1=-a_1-j+1=2$, and $b_1 \leq i_1 \leq -a_1$, so there is only one generator, 
namely $\frac{1}{\pi_1^2}$. Moreover, for $i\geq 2$ we have  $b_i=-a_j-j=-1$ and there are two possibilities 
for $-1 \leq i_\lambda \leq 0=-a_i$, namely $-1,0$. But only $\binom{0/n}{p-1}=0$, and we finally have that 
\[
H^1(V,\TO) \cong \langle \frac{1}{\pi_1^2} \rangle _k \times  \langle 1 \rangle _k \times \cdots \times  \langle 1 \rangle _k, 
\]
a space of dimension $\log_p |V|$. 

Let $d$ be the dimension of each simple direct summand of $H^1(V,\TO)$ considered as a $\Z/n\Z$-module. Of course 
$d$ equals the degree of the extenstion $\mathbb{F}_p(\zeta)/\mathbb{F}_p$, where $\zeta$ is a suitable primitive root of $1$.
For the matrix $\mathrm{diag}(\zeta^{c(1,i)},\ldots,\zeta^{c(d,i)})$ we have that 
\[
\mathrm{diag}(\zeta^{c(1,i)},\ldots,\zeta^{c(d,i)})=
\left\{
\begin{array}{ll}
\mathrm{diag}(\zeta^{2},\zeta,\ldots,\zeta) & \mbox{ if } i=1, \\
\zeta \cdot \mathbb{I}_d & \mbox{ if } i\geq 2.
\end{array}
\right.
\]
The characteristic polynomial in the first case is computed to be:
\[
f_1(x)=x^d + \sum_{\nu=1}^{d-1}  \zeta^{d-\nu} a_\nu x^\nu + a_0 \zeta^{1+d} . 
\]
By seting $x=\zeta y$ the above equation is written as
\[
\zeta^d\left(y^d + \sum_{\nu=0}^{d-1} a_\nu y^\nu \right) + (\zeta^{d+1}-\zeta^d)a_0.
\]
Therefore, for $x=1$, $y=\zeta$, so $f_1(1)=(\zeta^{d+1}-\zeta^d)a_0\neq 0$,
therefore $\delta(1)=0$.

In the second case, we observe that 
\[
f_i(x)=x^d + \sum_{\nu=0}^{d-1}  \zeta^{d-\nu}   a_\nu x^\nu.
\]
If we set $x=y/\zeta$, we obtain that $1$ is a simple root of  $f_i$, 
so $\delta(i)=1$ for $i\geq 2$.  

Thus, only the $s/d-1$ blocks $i\geq 2$ admit invariant elements and 
\[\dim_k H^1(V\rtimes \Z/n\Z,\TO)=s/d-1.\]

{\bf Comparison with the work of R. Pries.} 
Let us consider the curve:
\[
C: y^p-y =f(x), 
\]
where $f(x)$ is a polynomial of degree $j$, $(j,p)=1$. This, gives rise to a ramified cover of $\mathbb{P}^1(k)$ with $\infty$ as 
the unique ramification point. Moreover if all the monomial summands of the polynomial $f(x)$ have exponents congruent to $j\mod m,$
then the curve $C$ admits the group $G:=\Z/p\Z \rtimes \Z/m\Z$ as a subgroup of the group of automorphisms.
In \cite{Pries:02} R. Pries constructed a configuration  space of deformations of the above  curve and computed the dimension of this space. 
In what follows we will compute the dimension of $H^1(G,\TO)$ and we will compare with the result of Pries.

According to proposition  \ref{final-cyc-comp} the tangent space of the deformation space is generated as a $k$-vector space 
by the elements of the form $\frac{1}{x^i}$ where $b\leq i \leq j+1$ and 
\[
b=\left\{
\begin{array}{ll}
1 & \mbox{ if } p \mid -j-1\\
2 & \mbox{ if } p \nmid -j-1. 
\end{array}
\right.
\]
By lemma \ref{trivial-action} the $\Z/m\Z$-action on $\mathbb{F}_p$ is given by multiplication by $\zeta^j$ where $\zeta$ is an appropriate 
primitive $m$-th root of one. This gives us that $\zeta^{jp}=\zeta^j$, {\em i.e.} $jp\equiv j \mod m.$ 
If $d_i$ is the cocyle coresponding to the element $\frac{1}{x^i}$ then 
\[
d_i(\sigma \tau \sigma^{-1})^{\sigma^{-1}}=\zeta^j d_i(\tau)^{\sigma^{-1}}.
\]
But the element $\frac{1}{x^i}$ coresponds to the element $x^{j+1-i} \frac{d}{dx}$. The $\zeta^{-1}$-action is given by 
\[
x^{j+1-i}\frac{d}{dx} \mapsto \zeta^{i-j}x^{j+1-i}\frac{d}{dx},
\]
Therefore, the action of $\sigma$ on the cocycle corresponding to $\frac{1}{x^i}$ is given by 
\[
\frac{1}{x^i} \mapsto \zeta^i \frac{1}{x^i}.
\]
Thus,  $\dim_k H^1(G,\TO)=\dim_k H^1(\Z/p\Z,\TO)^{\Z/p\Z}$ is equal to 
\[
\#\left\{i: b\leq i \leq j+1,\; \binom{i/j}{p-1}=0,\; i\equiv 0 \mod m. \right\}
\] 
By (\ref{glob-contr}) we have 
\[ \dim_k H^1(Y,\pi_*^G (\T_X)) = 3 g_Y -3 + \sum_{\kappa=1}^r 
\lc \sum_{i=0}^{n_\kappa} \frac{(e_i^{(\kappa)}-1)}{e_0^{(\kappa)}} \rc,
\]
and by computation 
\[
\dim_k H^1(Y,\pi_*^G (\T_X)) = -3 + \lc 1+\frac{-1}{mp} + \frac{j(p-1)}{mp} \rc.
\]

On the other hand the configuration space constructed by R. Pries is of dimension 
\[
r:=\#\{e \in E_0:\forall \nu \in \mathbb{N}^+, p^\nu e \not \in E_0 \}
\]
where 
\[
E_0:=\{e: 1\leq e \leq j, e\equiv j \mod m.\}
\]
The  above dimensions look quite different but using {\em maple}\footnote{The maple code  used for this computation is available on
my web page {\tt http://eloris.samos.aegean.gr/preprints} for download}  we computed the following table:
\begin{center}
\begin{tabular}{l|l|l|l|l|l|l}
$p$ & $j$ & $m$ & $r$ & $\dim_k H^1(G,\TO)$ & $\dim_k  H^1(Y,\pi_*^G (\T_X))$ & $\dim_k D(k[\epsilon])$ \\
\hline
13 & 19 & 6  & 3 & 3 & 1 & 4 \\
13 & 35 & 6 &  5 & 4 & 9 & 13 \\ 
13 & 51 & 6 & 8 & 8 & 6 & 14 \\
13 & 36 & 3 & 12 & 11 & 10 & 21 \\
7 & 81 & 3 & 24 & 23 & 22 & 45 \\
7 & 90 & 3 & 26 & 26 & 24 & 50 \\
 \end{tabular}
\end{center}
We observe that the $r+a=\dim_k H^1(G,\TO)$, where $a=1,0$   and also the dimension of $H^1(Y,\pi_*^G (\T_X))$ is 
near the above two values.

Conclusions: in her paper R. Pries considered deformations, where the ramification point does not split to ramification 
points that are ramified with smaler jumps at their ramification filtrations.  By the difference of the formulas and by the 
fact that all infinitesimal deformations in 
 $H^1(Y,\pi_*^G (\T_X))$
 are unobstructed we obtain that the difference in the dimensions $r$ and 
 $\dim_k D(k[\epsilon])$ can be explained either as obstructed deformations or 
 as deformations splitting the branch points, see \ref{split-branch-points}.

 Let $P$ be a $p$-group. In \cite{Pries:04} R. Pries studied the  dimension of unobstructed deformations  of the 
 general case of $P\rtimes \Z/m\Z$, acting on a curve, without spliting the branch 
 points. A comparison of this dimension with the dimension computed by the tools 
 of this paper might be an interesting result,   although it is  a quite complicated task.  
  
\providecommand{\bysame}{\leavevmode\hbox to3em{\hrulefill}\thinspace}
\providecommand{\MR}{\relax\ifhmode\unskip\space\fi MR }
\providecommand{\MRhref}[2]{%
  \href{http://www.ams.org/mathscinet-getitem?mr=#1}{#2}
}
\providecommand{\href}[2]{#2}

\end{document}